\input amstex\documentstyle{amsppt}  
\pagewidth{12.5cm}\pageheight{19cm}\magnification\magstep1
\topmatter
\title Study of antiorbital complexes\endtitle
\author G. Lusztig\endauthor
\address{Department of Mathematics, M.I.T., Cambridge, MA 02139}\endaddress
\thanks{Supported in part by the National Science Foundation}\endthanks
\dedicatory{To Gregg Zuckerman on his 60th birthday}\enddedicatory
\endtopmatter   
\document    
\define\hcp{\hat{\cp}}

\define\tInd{\ti{\Ind}}
\define\tRes{\ti{\Res}}

\define\bco{\bar{\co}}

\define\uB{\un B}

\define\bXi{\bar\Xi}

\define\frl{\forall}
\define\pe{\perp}
\define\si{\sim}

\define\sqc{\sqcup}

\define\qua{\quad}

\define\hf{\hat f}

\define\bY{\bar Y}

\define\bsi{\bar\s}
\define\lb{\linebreak}

\define\op{\oplus}

\define\part{\partial}
\define\em{\emptyset}

\define\ra{\rangle}
\define\n{\notin}

\define\m{\mapsto}
\define\do{\dots}
\define\la{\langle}

\define\sub{\subset}    
\define\bxt{\boxtimes}
\define\T{\times}
\define\ti{\tilde}
\define\nl{\newline}
\redefine\i{^{-1}}
\define\fra{\frac}
\define\un{\underline}
\define\ov{\overline}
\define\ot{\otimes}
\define\bbq{\bar{\QQ}_l}

\define\Ad{\text{\rm Ad}}

\define\End{\text{\rm End}}

\define\Ind{\text{\rm Ind}}

\define\Res{\text{\rm Res}}

\define\tr{\text{\rm tr}}

\define\supp{\text{\rm supp}}

\define\a{\alpha}
\redefine\b{\beta}
\redefine\c{\chi}

\redefine\d{\delta}
\define\e{\epsilon}

\define\io{\iota}
\redefine\o{\omega}
\define\p{\pi}
\define\ph{\phi}
\define\ps{\psi}
\define\r{\rho}
\define\s{\sigma}
\redefine\t{\tau}
\define\th{\theta}
\define\k{\kappa}
\redefine\l{\lambda}
\define\z{\zeta}
\define\x{\xi}

\define\vt{\vartheta}

\redefine\G{\Gamma}

\redefine\L{\Lambda}
\define\Ph{\Phi}
\define\Ps{\Psi}

\define\kk{\bold k}

\define\vv{\bold v}

\define\FF{\bold F}

\define\LL{\bold L}

\define\NN{\bold N}

\define\QQ{\bold Q}

\define\SS{\bold S}

\define\VV{\bold V}

\define\ZZ{\bold Z}

\define\ca{\Cal A}
\define\cb{\Cal B}
\define\cc{\Cal C}
\define\cd{\Cal D}
\define\ce{\Cal E}
\define\cf{\Cal F}
\define\cg{\Cal G}

\define\ck{\Cal K}
\define\cl{\Cal L}

\define\co{\Cal O}
\define\cp{\Cal P}

\define\car{\Cal R}
\define\cs{\Cal S}

\define\cv{\Cal V}

\define\cx{\Cal X}
\define\cy{\Cal Y}

\define\fg{\frak g}
\define\fh{\frak h}

\define\fl{\frak l}

\define\fn{\frak n}

\define\fp{\frak p}

\define\fs{\frak s}

\define\fF{\frak F}

\define\fK{\frak K}

\define\fN{\frak N}

\define\fQ{\frak Q}

\define\fS{\frak S}

\define\tit{\ti t}

\define\tE{\ti E}

\define\tL{\ti L}

\define\tS{\ti S}

\define\tX{\ti X}

\define\tZ{\ti Z}

\define\sha{\sharp}

\define\uP{\un P}

\define\bS{\bar S}

\define\tce{\ti\ce}

\define\chE{\check E}

\define\DE{D}
\define\HE{H}
\define\LI{L2}
\define\LII{L3}
\define\LIII{L1}
\define\LIV{L4}
\define\LV{L5}
\define\VI{V}
\define\WE{W}

\head Introduction\endhead
\subhead 0.1\endsubhead
Let $U$ be an $N$-dimensional vector space over a finite field $\FF_q$ and let $U^*$ be the dual vector space. Let
$\k:U\T U^*@>>>\FF_q$ be the canonical pairing. We fix a prime number $l$ such that $l\ne0$ in $\FF_q$ and a
nontrivial character $\ps:\FF_q@>>>\bbq^*$. (We denote by $\bbq$ an algebraic closure of the field of $l$-adic 
numbers.) If $f:U@>>>\bbq$ is a function, the Fourier transform $\hf:U^*@>>>\bbq$ is defined by
$$\hf(y)=q^{-N/2}\sum_{x\in U}\ps(\k(x,y))f(x).$$

\subhead 0.2\endsubhead
Assume that $q$ is odd, $N$ is even $\ge4$ and that $U$ is endowed with a split nondegenerate symmetric bilinear 
form $(,):U\T U@>>>\FF_q$. This allows us to identify $U^*=U$. Let $SO(U)$ be the special orthogonal group of 
$(,)$. The following is the prototype of a problem that we are interested in.
 
(a) {\it Describe the space of all functions $f:U@>>>\bbq$ which are constant on each orbit of $SO(U)$ such that 
both $f$ and $\hf$ vanish on the complement of $\{x\in U;(x,x)=0\}$.}
\nl
It turns out that the space (a) consists of all scalar multiples of a single function namely:
$$\align&f(x)=0\text{ if } (x,x)\ne0,\qua f(x)=1\text{ if }(x,x)=0,x\ne0,\\&
f(x)=1+q^{(N-2)/2}\text{ if }x=0.\endalign$$
The proof is an easy exercise (the fact that $\hf=f$ is shown in \cite{\LIV,\S11}; see also the last paragraph of
3.5). 

Let $\kk$ be an algebraic closure of $\FF_q$. The function $f$ can be interpreted as the characteristic function 
of the intersection cohomology complex of the quadric $(x,x)=0$ in $\kk\ot U$. This quadric is singular at $0$ and
this accounts for the term $q^{(N-2)/2}$ in the formula for $f(0)$.

\subhead 0.3\endsubhead
Let $G$ be a connected reductive algebraic group over $\kk$ with a given semisimple automorphism $\vt:G@>>>G$ with
fixed point set $G^\vt$. Let $K$ be the identity component of $G^\vt$. We assume that the characteristic of $\kk$
is sufficiently large. Now $\vt$ induces a (semisimple) automorphism of the Lie algebra $\fg$ of $G$. For any 
$\z\in\kk^*$ let $\fg_\z$ be the $\z$-eigenspace of this automorphism. Note that $K$ acts on $\fg_\z$ by the 
adjoint action. Let $\fg_\z^{nil}$ be the variety of elements in $\fg_\z$ which are nilpotent in $\fg$. We assume
that we are given an $\FF_q$-structure on $G$ such that $\vt:G@>>>G$ is defined over $\FF_q$, such that all 
eigenvalues of $\vt:\fg@>>>\fg$ are in $\FF_q^*$ and such that there exists a Borel subgroup of $G$ defined over 
$\FF_q$ which is $\vt$-stable. Then $K$, $\fg$ and $\fg_\z$ (for any $\z$) inherit an $\FF_q$-structure. Let 
$\k:\fg\T\fg@>>>\kk$ be a fixed $G$-invariant nondegenerate symmetric bilinear form; we assume that it is also
$\vt$-invariant and defined over $\FF_q$.
It induces a nondegenerate bilinear pairing $\k_\z:\fg_\z\T\fg_{\z\i}@>>>\kk$ (for any $\z$) which is again 
defined over $\FF_q$. This allows us to identify $\fg_{\z\i}$ with the dual space of $\fg_\z$ (compatibly with the
$\FF_q$-structures and $K$-actions). We fix $\z\in\kk^*$. The following is a generalization of the problem 0.2(a).

(a) {\it Describe the space of all functions $f:\fg_\z(\FF_q)@>>>\bbq$ (which are constant on each orbit of 
$K(\FF_q)$) such that $f$ vanishes on any non-nilpotent element of $\fg_\z(\FF_q)$ and $\hf$ vanishes on any 
non-nilpotent element of $\fg_{\z\i}(\FF_q)$. In particular, when is this space nonzero?}
\nl
(The example in 0.2 arises in the case where $G=SO_{2n+1}(\kk)$, $K=SO_{2n}(\kk)$, $\z=-1$.)

In the case where $\z=1$ the solution of (a) can be found in \cite{\LIII}. It turns out that in this case
functions as in (a) exist very rarely and they are related to "cuspidal character sheaves". For example if 
$G=Sp_{2n}(\kk)$, $\vt=1,\z=1$, there is up to scalar at most one nonzero $f$ as in (a); it exists if and only if
$n=i(i+1)/2$ for some $i\in\NN$. On the other hand, when $\vt$ is an inner automorphism and $\z\ne1$, we will show
(see 3.4(a)) that the space (a) is always nonzero. 

\subhead 0.4\endsubhead
Assume now that in 0.3 we take $G=GL(V)$ where $V$ is a finite dimensional $\kk$-vector space with a fixed
$\FF_q$-structure and a fixed grading $V=\op_{i\in\ZZ/m}V_i$ such that each $V_i$ is defined over $\FF_q$. Here 
$m$ is a fixed integer $\ge2$ such that $m\ne0$ in $\kk$. Let $a\in\kk^*$ be such that $a^m=1$ and 
$1,a,\do,a^{m-1}$ are distinct. Define $g_0\in G$ by $g_0(v)=a^iv$ for $v\in V_i,i\in\ZZ/m$. Define $\vt:G@>>>G$ 
as conjugation by $g_0$. Then $\fg_a$ can be viewed as the space of representations of a cyclic quiver with $m$ 
vertices of fixed multidegree. In this case we will describe completely the space 0.3(a) by exhibiting an explicit
basis for it. (See 2.14(a)). This basis is closely connected with the theory of canonical bases \cite{\LI} applied
to affine $SL_m$.

\subhead 0.5\endsubhead
As suggested by the example in 0.2, it is useful to study the problem 0.3(a) by passing to $\kk$ and using 
geometric methods, namely using perverse sheaves instead of functions and using Deligne-Fourier (D-F) transform 
(see 0.7) instead of Fourier transform.

Let $E$ be a finite dimensional $\kk$-vector space and let $\chE$ be its dual space. Let $K$ be a connected linear
algebraic group with a given homomorphism of algebraic groups $K@>>>GL(E)$. Then $K$ acts on $E$ and (by duality) 
on $\chE$. A complex (of $\bbq$-sheaves) on $E$ is said to be {\it orbital} if it is a simple $K$-equivariant 
perverse sheaf on $E$ supported by the closure of a single $K$-orbit. A complex (of $\bbq$-sheaves) on $\chE$ is
said to be {\it antiorbital} if it is isomorphic to the D-F transform of an orbital complex on $E$. (Then it is a
simple $K$-equivariant perverse sheaf on $\chE$.) A complex (of $\bbq$-sheaves) on $E$ is said to be {\it 
biorbital} if it is orbital for the $K$-action on $E$ and its D-F transform is orbital for the $K$-action on
$\chE$. 

\subhead 0.6\endsubhead
Assume now that $E=\fg_\z,K$ are as in 0.3, ($\z\in\kk^*$). Some general results on $K$-orbits on $\fg_\z$ can be
deduced from Vinberg's work \cite{\VI}. In the case where $\z=1$ an explicit description of the antiorbital 
complexes on $\fg_{\z\i}$ is given in \cite{\LIII}. 

In the case where $\fg_\z$ arises
as in 0.4 from a cyclic quiver an explicit description of the antiorbital complexes on $\fg_{\z\i}$ is given in 
Theorem 2.6. By "explicit" we mean that for each antiorbital complex we describe its support and the associated 
local system on an open dense smooth part of the support. 

We also describe explicitly the collection of antiorbital complexes on $\chE$ in a case where $E,K$ does not come
from the setup in 0.3 (see 3.9). In this example, $E$ is the coadjoint representation of a unipotent group $K$ and
there are infinitely many subvarieties of the adjoint representation $\chE$ which appear as supports of 
antiorbital complexes.

Assume again that $E=\fg_\z,K,\vt$ are as in 0.3 with $\z\in\kk^*$, so that $\chE=\fg_{\z\i}$. Let $\fQ_\z$ be the
collection of simple $K$-equivariant perverse sheaves $\ck$ on $\fg_\z$ such that $\supp(\ck)\sub\fg_\z^{nil}$ and
$\supp(\fF(\ck))\sub\fg_{\z\i}^{nil}$. Note that $\fQ_\z$ contains only finitely many objects up to isomorphism.
The following is a geometric analogue of problem 0.3(a).

(a) {\it Describe explicitly the objects of $\fQ_\z$ (by describing the corresponding nilpotent $K$-orbits and 
associated local systems). In particular when is $\fQ_\z$ nonempty?}
\nl
In the case where $\z=1$ the solution of this problem can be found in \cite{\LIII}. It turns out that in this case
$\fQ_\z$ is almost always empty. For example if $G=Sp_{2n}(\kk)$, $\vt=1,\z=1$, there is up to isomorphism at most
one object of $\fQ_\z$; it exists if and only if $n=i(i+1)/2$ for some $i\in\NN$. At the other extreme, in the 
case where $\z$ has large order in $\kk^*$ (say $>\dim\fg$), $\fg_{\z\i}$ and $\fg_\z$ consist of nilpotent 
elements and $\fQ_\z$ consists of all orbital complexes on $\fg_{\z\i}$ (they are automatically antiorbital). If 
we assume only that $\z\ne1$ (and that $\vt$ is an inner automorphism), we can show that $\fQ_\z\ne\em$. (See 
3.3.) In the case which arises as in 0.4 from a cyclic quiver, we determine $\fQ_\z$ explicitly. (See 2.13). In 
this case $\fQ_\z$ is exactly the collection of simple perverse sheaves on $E$ defined by the theory of canonical
bases for affine $SL_m$. (This gives a new characterization of that canonical basis.) In any case, the set (of 
isomorphism classes in) $\fQ_\z$ can be viewed as a basis of a finite dimensional $\QQ(v)$-vector space $\VV_\z$ 
and we have the operations of "induction" and "restriction" (see the Appendix) relating $\VV_\z$ and the analogous
vector space for a $\vt$-stable Levi subgroup analogous to multiplication and comultiplication in a quantum group;
we also have canonically $\VV_\z=\VV_{\z\i}$ (compatibly with the bases) and this agrees with the operations of 
induction and restriction. Thus, $\fQ_\z$ behaves in many respects like a canonical basis. (This point of view has
been already emphasized in \cite{\LV} which corresponds to the case where $\z$ has large order in $\kk^*$.)

\subhead 0.7. Notation\endsubhead
If $X$ is an algebraic variety let $\cd(X)$ be the bounded derived category of constructible $\bbq$-sheaves on $X$
The objects of $\cd(X)$ are said to be complexes. If $\ck\in\cd(X)$ we write $\ck[?]$ instead of "$\ck$ with some
shift". We say "local system" instead of "$\bbq$-local system". If $\cl$
is a local system on $X$ and $x\in X$ we denote by $\cl_x$ the stalk of $\cl$ at $x$. If $Y$ is a locally closed 
smooth irreducible subvariety of $X$ and $\cl$ is a local system on $Y$, then the intersection cohomology complex
$IC(\bY,\cl)$ on the closure $\bY$ of $Y$ is well defined. Extending this by $0$ on $X-\bY$ we obtain a complex 
$IC(\bY,\cl)_X$ in $\cd(X)$. Now $\ps$ (see 0.1) defines an Artin-Schreier local system of rank $1$ on $\kk$; its
inverse image under any morphism $\ph:X@>>>\kk$ is a local system $\cl^\ph$ of rank $1$ on $X$.

Let $E,\chE$ be as in 0.5. Let $\la,\ra:E\T\chE@>>>\kk$ be the obvious pairing. Let $\fF:\cd(E)@>>>\cd(\chE)$ be 
the Deligne-Fourier (D-F) transform 
$$\ck\m\fF(\ck)=\d'_!(\d^*(\ck)\ot\cl^{\la,\ra})[\dim E];$$
here $\d:E\T\chE@>>>E$, $\d':E\T\chE@>>>\chE$ are the two projections. 

For $n\in\NN$, $\fS_n$ denotes the symmetric group in $n$ letters.

\head Contents\endhead
1. Cyclic quivers and orbital complexes.

2. Cyclic quivers and antiorbital complexes.

3. Further examples of antiorbital complexes.

Appendix. Induction, restriction.

\head 1. Cyclic quivers and orbital complexes\endhead
\subhead 1.1\endsubhead
In this and the next section we fix an integer $m\ge1$ such that $m\ne0$ in $\kk$. We set $I=\ZZ/m$. We fix 
$\e=\pm1$ in $I$.

Let $\cc$ be the category of finite dimensional $\kk$-vector spaces $V$ with a given $I$-grading that is, a direct
sum decomposition $V=\op_{i\in I}V_i$ indexed by $I$ (the morphisms are linear maps compatible with the grading).
For $V\in\cc$ we define $|V|=(|V|_i)\in\NN^I$ by $|V|_i=\dim(V_i)$. Define $\vv\in\cc$ by $\vv_i=\kk$ for all 
$i\in I$. For $V\in\cc$ we set 
$$E^\e_V=\{T\in\End(V);TV_i\sub V_{i+\e}\qua\frl i\}.$$
Let $\cc^\e$ be the category whose objects are the pairs $(V,T)$ where $V\in\cc$ and $T\in E^\e_V$ (the morphisms
are linear maps compatible with the grading and which commute with the given endomorphisms). For $V\in\cc$ let 
$$G_V=\{A\in GL(V),AV_i=V_i\qua\frl i\in I\},$$
a closed subgroup of $GL(V)$ isomorphic to $\prod_{i\in I}GL(V_i)$.
Note that $G_V$ acts naturally by conjugation on $E^\e_V$. Clearly the $G_V$-orbits on $E^\e_V$ are in natural 
bijection with the isomorphism classes of objects $(V',T)\in\cc^\e$ such that $V'\cong V$ in $\cc$. 

For any $T\in E^\e_V$ and any $\l\in\kk$ let $V_{T,\l}$ be the generalized $\l$-eigenspace of $T^m:V@>>>V$. Note 
that $V_{T,\l}$ is compatible with the grading of $V$ hence is itself an object of $\cc$. Moreover $V_{T,\l}$ is 
$T$-stable. We denote the restriction of $T$ to $V_{T,\l}$ by ${}^\l T$. We have $(V_{T,\l},{}^\l T)\in\cc^\e$ and
$(V,T)=\op_{\l\in\kk}(V_{T,\l},{}^\l T)$ in $\cc^\e$. We set 
$$\cs_T=\{\l\in\kk^*;V_{T,\l}\ne0\}.$$
For any $\l\in\kk$ let $E^{\e,\l}_V=\{T\in E^\e_V;V_{T,\l}=V\}$ (a closed $G_V$-stable subset of $E^\e_V$). 
We say that $E^{\e,0}_V$ is the {\it nilpotent variety} in $E^\e_V$.

If $V\cong\vv$ we set $E^{-\e,*}_V=E^{-\e}_V-E^{-\e,0}_V=\cup_{\l\in\kk^*}E^{-\e,\l}_V$ (an open dense subset of 
$E^{-\e}_V$); for $\l\in\kk^*$ we define $\a_\l:E^{-\e,*}_V@>>>\kk^*$ by $\a_\l(T)=\l\l'$ where 
$T\in E^{-\e,\l'}_V$.

We describe some objects of $\cc^\e$.

Let $\ca=\{(a,b)\in\ZZ\T\ZZ;a\le b\}$. Now $m\ZZ$ acts freely on $\ca$ by $c:(a,b)\m(a+c,b+c)$; let $\ca_m$ be 
the space of orbits. We write $\ov{a,b}$ for the orbit of $(a,b)\in\ca$.

Let $(a,b)\in\ca$. Let $V=V_{a,b}$ be the $\kk$-vector space with basis $e_a,e_{a+1},\do,e_b$ viewed as an object
of $\cc$ in which $e_j\in V_j$ for all $j=a,a+1,\do,b$. (The index $j$ in $V_j$ is viewed as an integer $\mod m$.)
Define $T_1\in E^1_V$ by $e_a\m e_{a+1}\m e_{a+2}\m\do\m e_b\m0$. Define $T_{-1}\in E^{-1}_V$ by 
$e_b\m e_{b-1}\m e_{b-2}\m\do\m e_a\m0$. We have $(V_{a,b},T_\e)\in\cc^\e$. Moreover the isomorphism class of
$(V_{a,b},T_\e)$ in $\cc^\e$ depends only on the $m\ZZ$-orbit of $(a,b)$. Hence we can use the notation
$(V_\b,T_\e)$ for this isomorphism class where $\b=\ov{a,b}\in\ca_m$.

Let $\l\in\kk^*$ and let $n\in\ZZ_{>0}$. Choose $\l'\in\kk^*$ such that $\l'{}^m=\l$. Let $U=U(n)$ be the
$\kk$-vector space with basis $e_{a,i}$ $(a\in[1,n],i\in I)$. For $i\in I$ let $U_i$ be the subspace of $U$ 
spanned by $e_{a,i}(a\in[1,n])$. These subspaces form an $I$-grading of $U$. Define $T_1\in\End(U)$ by 
$e_{a,i}\m e_{a,i+\e}$ $(a\in[1,n],i\in I)$. Define $T_2\in\End(U)$ by $e_{a,i}\m e_{a+1,i+\e}$ 
$(a\in[1,n-1],i\in I)$, $e_{n,i}\m0$ $(i\in I)$. We have $T_1,T_2\in E^\e_{U}$, $T_1T_2=T_2T_1$. We have 
$T_1^m=1$ hence $T_1$ is semisimple. Clearly, $T_2$ is nilpotent hence 
$T_\e(\l):=\l'T_1+T_2\in\End(U)$ is such that $T_\e(\l)^m-\l:U@>>>U$ is nilpotent. We have 
$(U(n),T_\e(\l))\in\cc^\e$. We show that the isomorphism class of $(U,T_\e(\l))$ is independent of the choice of 
$\l'$. Let $\l''\in\kk^*$ be such that $\l''{}^m=\l$. We must show that there exists $R\in G_U$ such that 
$R(\l'T_1+T_2)=(\l''T_1+T_2)R$. Let $z=\l''/\l'$. Define $R\in G_U$ by $e_{a,i}\m z^{-a}z^{\e i}e_{a,i}$ 
$(a\in[1,n],i\in I)$. (Here $z^{\e i}$ makes sense since $z^m=1$.) We have $RT_1=zT_1$, $RT_2=T_2R$ and our claim
follows.

Now $(V_\b,T_\e)$ and $(U(n),T_\e(\l))$ are indecomposable objects of $\cc^\e$; it is easy to see that,
conversely, any indecomposable object of $\cc^\e$ is isomorphic to $(V_\b,T_\e)$ (for a well defined $\b\in\ca_m$)
or to $(U(n),T_\e(\l))$ (for a well defined $n\ge1$, $\l\in\kk^*$).

\subhead 1.2\endsubhead
Define $\bold1\in\NN^I$ by $\bold1_i=1$ for all $i\in I$. Let $\car$ be the set of maps $\r:\ZZ_{>0}@>>>\NN$ such
that $\r(n)=0$ for all but finitely many $n$. For $\r\in\car$ we set $\un\r=\sum_{n\in\ZZ_{>0}}\r(n)n\in\NN$. For
any $t\in\NN$ let $\car_t=\{\r\in\car;\un\r=t\}$.

Let $\cp$ be the set of maps $\s:\ca_m@>>>\NN$ such that $\s(\b)=0$ for all but finitely many $\b\in\ca_m$. Let
$$\cp^{ap}=\{\s\in\cp;\prod_{\ov{a,b}\in\ca_m;b-a=n}\s(\ov{a,b})=0\qua\frl n\in\NN\}.$$ 
(This is a finite product.) The elements of $\cp^{ap}$ are said to be {\it aperiodic}. For $\s\in\cp$ we set 
$|\s|=\sum_{\b\in\ca_m}\s(\b)|V_\b|\in\NN^I$. For any $\nu\in\NN^I$ we set $\cp_\nu=\{\s\in\cp;|\s|=\nu\}$. Let 
$\cp^{ap}_\nu=\{\s\in\cp^{ap};|\s|=\nu\}$.

Let $\tZ$ be the set of all collections $\th=(\th_\l)_{\l\in\kk}$ where $\th_\l\in\car$ for $\l\in\kk^*$, 
$\th_0\in\cp$ are such that $\th_\l=0$ for all but finitely many $\l$. For $\th\in\tZ$ we set 
$|\th|=\sum_{\l\in\kk^*}\un{\th_\l}\bold1+|\th_0|\in\NN^I$. For $\nu\in\NN^I$ let 
$\tZ_\nu=\{\th\in\tZ;|\th|=\nu\}$.

For any $\r\in\car$ let $\Pi_\r$ be the corresponding irreducible representation of $\fS_{\un\r}$ over $\bbq$. (In
particular if $\r(1)=N$ and $\r(n)=0$ for $n>1$ then $\Pi_\r$ is the unit representation of $\fS_N$.) Let 
$N_\r=\dim\Pi_\r$.

\subhead 1.3\endsubhead
Let $V\in\cc$. From the results in 1.1 we see that the $G_V$-orbits on $E^\e_V$ are naturally indexed by the set 
$\tZ_{|V|}$: the $G_V$-orbit $E^\e_{V,\th}$ corresponding to $\th\in\tZ_{|V|}$ consists of all $T\in E^\e_V$ such
that
$$(V,T)\cong\op_{\l\in\kk^*,n\ge1}(U(n),T_\e(\l))^{\op\th_\l(n)}\op\op_{\b\in\ca_m}(V_\b,T_\e)^{\op\th_0(\b)}$$
in $\cc^\e$. For any $\th\in\tZ_{|V|}$, the complex
$$IC(\ov{E^\e_{V,\th}},\bbq)_{E^\e_V}\tag a$$
is well defined. It is up (to shift) an orbital complex on $E^\e_V$ with its natural $G_V$-action. Since the only
$G_V$-equivariant irreducible local system on $E^\e_{V,\th}$ is $\bbq$ we see that any orbital complex on $E^\e_V$
is (up to shift) of the form (a) for a well defined $\th$.

Let $\s\in\cp_{|V|}$. Define $\bsi\in\tZ_{|V|}$ by $\bsi_0=\s$, $\bsi_\l=0$ for $\l\in\kk^*$. Note that the
$G_V$-orbits contained in the nilpotent variety $E^{\e,0}_V$ are exactly the subsets $E^\e_{V,\bsi}$ for various
$\s\in\cp_{|V|}$. 

\subhead 1.4\endsubhead
Let $V\in\cc$ be such that $V\cong\vv^{\op s}$ for some $s\in\NN$. We fix $\l\in\kk^*$. Let 
$D=\{S\in E^\e_V;S^m=\l\}$. For $S\in D$ let $\dot S=\{N\in E^{\e,0}_V;NS=NS\}$. Clearly, the map 
$$\{(S,N);S\in D,N\in\dot S\}@>g>>E^{\e,\l}_V,\qua (S,N)\m S+N$$
is a bijection. Let $\fN=\{f\in\End(V_0);f\text{ nilpotent}\}$. For any $S\in D$ we define a bijection 
$\fN@>\si>>\dot S$ by $\k\m N=N_{S,\k}$ where $N(S^hv)=S^{h+1}\k(v)$ for $v\in V_0$, $h\in[0,m-1]$; we then have 
also $N(S^hv)=S^{h+1}\k(v)$ for $v\in V_0$, $h\in\NN$. The inverse bijection is $N\m\k$ where $\k(v_0)=S\i N(v_0)$
for $v_0\in V_0$. We see that we have a bijection $u:D\T\fN@>\si>>E^{\e,\l}_V$ given by $(S,\k)\m S+N_{S,\k}$. Now
let $\SS$ be the set of all $(T,f)$ where $f=(0=V^0\sub V^1\sub V^2\sub\do\sub V^s=V)$ is a sequence of subobjects
of $V$ such that $V^j/V^{j-1}\cong\vv$ for $j\in[1,s]$ and $T\in E^\e_V$ is such that $T(V^j)\sub V^j$ for 
$j\in[1,s]$ and the element of $E^\e_{V^j/V^{j-1}}$ induced by $T$ is in $E^{\e,\l}_{V^j/V^{j-1}}$. This is a
smooth irreducible variety. Let $\SS_0$ be the set of all $(\k,f_0)$ where 
$f_0=(0=\cv^0\sub\cv^1\sub\cv^2\sub\do\sub\cv^s=V_0)$ is a sequence of subspaces of $V_0$ such that 
$\dim(\cv^j/\cv^{j-1})=1$ for $j\in[1,s]$ and $\k\in\fN$ is such that $\k(\cv^j)\sub\cv^j$ for $j\in[1,s]$. We 
define an isomorphism $u':D\T\SS_0@>\si>>\SS$ by $(S,(\k,f_0))\m(T,f)$ where $T=S+N_{S,\k}$ and $f=(V^j)$ with 
$V^j_h=S^h(\cv^j)$ for $h=0,1,\do,m-1$. We have a commutative diagram
$$\CD
\SS@<u'<< D\T\SS_0@>a'>>\SS_0 \\
@Vg_1VV     @Vg_2VV  @Vg_3VV \\
E^{\e,\l}_V@<u<<D\T\fN@>a>>\fN
\endCD$$
where $g_1(T,f)=T$, $g_2(S,(\k,f_0))=(S,\k)$, $g_3(\k,f_0)=\k$ and $a,a'$ are the obvious projections. It is well
known that $g_{3!}\bbq\cong\op_{\r\in\car_s}P_\r^{\op N_\r}[d]$ where $P_\r$ are simple mutually nonisomorphic
perverse sheaves on $\fN$ and $d\in\ZZ$. Using the fact that the right square in the diagram above is cartesian 
and $a$ is smooth with connected fibres we deduce that $g_{2!}\bbq\cong\op_{\r\in\car_s}P'_\r{}^{\op N_\r}[d']$ 
where $P'_\r$ are simple mutually nonisomorphic perverse sheaves on $S\T\fN$ and $d'\in\ZZ$. Using the fact that 
$u,u'$ are isomorphisms we deduce that $g_{1!}\bbq\cong\op_{\r\in\car_s}P''_\r{}^{\op N_\r}[d']$ where $P''_\r$ 
are simple mutually nonisomorphic perverse sheaves on $E^{\e,\l}_V$ and $d'$ is as above.

\subhead 1.5\endsubhead
Assume that $V\in\cc$, $U^1,U^2,\do,U^s\in\cc$ are such that $U^1\op U^2\op\do\op U^s\cong V$. We have a diagram
$$E_{U^1}^\e\T E_{U^2}^\e\T\do\T E_{U^s}^\e@<p_1<<E'@>p_2>>E''@>p_3>>>E^\e_V$$
with the following notation. 

$E''$ is the set of all pairs $(T,f)$ where $f$ is a sequence $(0=V^0\sub V^1\sub V^2\sub\do\sub V^s=V)$ of 
subobjects of $V$ such that $V^j/V^{j-1}\cong U^j$ for $j\in[1,s]$ and $T\in E^\e_V$ is such that $TV^j\sub V^j$ 
for all $j\in[1,s]$. $E'$ is the set of all triples $(T,f,\ph)$ where $(T,f)\in E''$ (with $f=(V^j)$) and 
$\ph=(\ph_j)$ is a collection of isomorphisms $\ph_j:V^j/V^{j-1}@>\si>>U^j$ for $j\in[1,s]$. We have 
$p_1(T,f,\ph)=(T_j)$ where $T_j\in E^\e_{U^j}$ is obtained by transporting the element of $E^\e_{V^j/V^{j-1}}$ 
induced by $T$ via $\ph_j$ for $j\in[1,s]$. We have $p_2(T,f,\ph)=(T,f)$, $p_3(T,f)=T$. Note that $p_3$ is a 
proper morphism.

Let $\ck_j$ be a $G_{U^j}$-equivariant perverse sheaf (up to shift) on $E^\e_{U^j}$ ($j\in[1,s])$. Then 
$\bxt_j\ck_j$ is a $\prod_jG_{U^j}$-equivariant perverse sheaf (up to shift) on $\prod_jE^\e_{U^j}$. Since $p_1$ 
is a smooth morphism with connected fibres, $\ck':=p_1^*(\bxt_j\ck_j)$ is a $\prod_jG_{U^j}$-equivariant perverse
sheaf (up to shift) on $E'$. Since $p_2$ is a principal $\prod_jG_{U^j}$-bundle there is a well defined perverse 
sheaf (up to shift) $\ck''$ on $E''$ such that $p_2^*\ck''=p_1^*\ck'$. We set 
$\Ind(\bxt_j\ck_j)=p_{3!}\ck''\in\cd(E^\e_V)$.

\subhead 1.6\endsubhead
We preserve the setup of 1.5. We fix $\l\in\kk^*$. Assume that $U^j\cong\vv$ for $j\in[1,s]$. For $j\in[1,s]$ let 
$\ck_j=IC(E^{\e,\l}_{U^j},\bbq)_{E^\e_{U^j}}$. Let $\SS$ be the set of all $(T,f)\in E''$ such that the following 
holds: $f=(V^j)$ is such that for $j\in[1,s]$, $V^j/V^{j-1}\cong U^j$ and the element of $E^\e_{V^j/V^{j-1}}$ 
induced by $T$ is in $E^{\e,\l}_{V^j/V^{j-1}}$. This is a closed smooth irreducible subset of $E''$. From the 
definitions we have $\ck''=IC(\SS,\bbq)_{E''}$. We have $p_3(\SS)\sub E^{\e,\l}_V$. From 1.4 we see that 
$$\Ind(\bxt_j\ck_j)\cong\op_{\r\in\car_s}P''_\r{}^{\op N_\r}[d']$$
where $P''_\r$ are simple mutually nonisomorphic perverse sheaves on $E^\e_V$ with support on $E^{\e,\l}_V$ and 
$d'\in\ZZ$. The perverse sheaves $P''_\r$ are $G_V$-equivariant. Since the number of $G_V$-orbits in $E^{\e,\l}_V$
is equal to $\sha(\car_s)$ we see that the $P''_\r$ are precisely the orbital complexes on $E^\e_V$ with support 
contained in $E^{\e,\l}_V$.

\subhead 1.7\endsubhead
We preserve the setup of 1.5. Assume that $s=t+1$ where $t\in\NN$. Assume that $U^j\cong\vv$ for $j\in[1,t]$ and 
$\s\in\cp^{ap}$ is such that $|U^s|=|\s|$. For $j\in[1,t]$ let $\ck_j=IC(\{0\},\bbq)_{E^\e_{U^j}}$. Let 
$\ck_s=IC(\ov{E^\e_{U^s,\bsi}},\bbq)_{E^\e_{U^s}}$. Let $S''$ be the set of all $(T,f)\in E''$ such that the 
following holds: $f=(V^j)$ is such that $V^j/V^{j-1}\cong U^j$ for $j\in[1,s]$ and the element of 
$E^\e_{V^j/V^{j-1}}$ induced by $T$ is $0$ if $j\in[1,t]$ and is in $E^\e_{V^s/V^{s-1},\bsi}$ if $j=s$. This is a
locally closed smooth irreducible subset of $E''$. From the definitions we have $\ck''=IC(\bS'',\bbq)_{E''}$. We
have $p_3(S'')\sub E^\e_{V,0}$. Since $E^\e_{V,0}$ is closed in $E^\e_V$ we see that $p_3(\bS'')\sub E^{\e,0}_V$.
Thus $\supp(\Ind(\bxt_j\ck_j))\sub E^{\e,0}_V$. Moreover, since $p_3$ is proper, $\Ind(\bxt_j\ck_j)$ is isomorphic
to $\op_{e\in[1,M]}P_e[d_e]$ where $P_e$ are simple $G_V$-equivariant perverse sheaves on $E^\e_V$ with 
$\supp(P_e)\sub E^{\e,0}_V$ and $d_e\in\ZZ$ (we use the decomposition theorem). Let $\cy_{t,\s}$ be the set of 
isomorphism classes of simple perverse sheaves on $E^\e_V$ that are isomorphic to $P_e$ for some $e\in[1,M]$.

\subhead 1.8\endsubhead
We preserve the setup of 1.5. Assume that $s=t+1$, $t\in\NN$. Let $\th\in\tZ_{|V|}$. Assume that for 
$j\in[1,t]$ we have $|U^j|=\un{\th_{\l_j}}\bold1$ where $\l_1,\l_2,\do,\l_t$ are distinct elements of $\kk^*$ and
$|U^s|=|\th_{\l_s}|$ (we set $\l_s=0$). Note that $\th_\l=0$ for any $\l\n\{\l_1,\do,\l_s\}$. For $j\in[1,s]$ 
define $\th^j\in\tZ_{|U^j|}$ by $\th^j_{\l_j}=\th_{\l_j}$, $\th^j_\l=0$ if $\l\ne\l^j$. Let 
$S_j=E^\e_{U^j,\th^j}$, a $G_{U^j}$-orbit in $E^\e_{U^j}$; let $\ck_j=IC(\bS_j,\bbq)_{E^\e_{U^j}}$. Let $S''$ be 
the set of all $(T,f)\in E''$ such that the following holds: $f=(V^j)$ is such that $V^j/V^{j-1}\cong U^j$ for 
$j\in[1,s]$ and $T\in E^\e_V$ is such that for $j\in[1,s]$, the element of $E^\e_{V^j/V^{j-1}}$ induced by $T$ 
belongs to $E^\e_{V^j/V^{j-1},\th^j}$. We have $\ck''=IC(\bS'',\bbq)_{E''}$. Let $S$ be the set of all 
$T\in E^\e_V$ such that for $j\in[1,s]$ we have $V_{T,\l_j}\cong U^j$ and the element of $E^\e_{V_{T,\l_j}}$ 
induced by $T$ belongs to $E^\e_{V_{T,\l_j},\th^j}$. Note that $S=E^\e_{V,\th}$. For any $T\in\bS$ (closure of 
$S$) the fibre $p_3\i(T)$ consists of exactly one element, namely $(T,f)$ where $f=(V^j)$ is given by 
$V^1=V_{T,\l_1}$, $V^2=V_{T,\l_1}\op V_{T,\l_2}$, etc. More precisely, $p_3$ defines an isomorphism 
$\bS''@>\si>>\bS$. We see that $\Ind(\bxt_j\ck_j)=IC(\ov{E^\e_{V,\th}},\bbq)_{E^\e_V}$. 

\head 2. Cyclic quivers and antiorbital complexes\endhead
\subhead 2.1\endsubhead
We preserve the setup in 1.1. If $V\in\cc$, we have a perfect bilinear pairing $E^\e_V\T E^{-\e}_V@>>>\kk$ given 
by $(T,T')=\tr(TT',V)$; it is compatible with the $G_V$-action. We use it to identify $E^{-\e}_V$ with the dual of
$E^\e_V$. Hence the D-F transform $\fF:\cd(E^\e_V)@>>>\cd(E^{-\e}_V)$ is well defined. In particular the notion of
antiorbital complex on $E^{-\e}_V$ is well defined; it is a complex of the form $\fF(\ck)$ where $\ck$ is an 
orbital complex on $E^\e_V$ (with respect to the $G_V$-action).

\subhead 2.2\endsubhead
Define $h:(\kk^*)^m@>>>\kk$ and $h':(\kk^*)^m@>>>\kk^*$ by 
$$h(x_1,x_2,\do,x_m)=x_1+x_2+\do+x_m,\qua h'(x_1,x_2,\do,x_m)=x_1x_2\do x_m.$$
Then the local system $\cl^h$ on $(\kk^*)^m$ is defined (see 0.7). According to Deligne 
\cite{\DE, Thm.7.8, p.221}, the complex $\fK^m:=h'_!\cl^h[m-1]\in\cd(\kk^*)$ (a sheaf theory version of a family 
of generalized Kloosterman sums) is a local system of rank $m$ on $\kk^*$. (The rank $1$ local system $\fK^1$ is
of Artin-Schreier type. The rank $2$ local system $\fK^2$ is implicit in Weil's paper \cite{\WE}.)

Assume now that $V\in\cc$, $V\cong\vv$. Let $\l\in\kk^*$. Let $\ck=IC(E^{\e,\l}_V,\bbq)_{E^\e_V}$. Note that 
$\a_\l^*(\fK^m)$ is a local system of rank $m$ on $E^{-\e,*}_V$ ($\a_\l$ as in 1.1). Let 
$\ck'=IC(E^{-\e}_V,\a_\l^*(\fK^m))_{E^{-\e}_V}$. We show:
$$\ck'[m]\cong\fF(\ck)[m-1].\tag a$$
Since $\fF(\ck)[m-1]$ is a simple perverse sheaf on $E^{-\e}_V$ it is enough to show that 
$$\fF(\ck)[m-1]|_{E^{-\e,*}_V}\cong\ck'[m]|_{E^{-\e,*}_V}$$
in $\cd(E^{-\e,*}_V)$ or equivalently that
$$\d'_{1!}(\d^*_1\bbq\ot\cl^{(,)})[m][m-1]=\a_\l^*(\fK^m)[m]$$
where $\d_1:E^{\e,\l}_V\T E^{-\e,*}_V@>>>E^{\e,\l}_V$, $\d'_1:E^{\e,\l}_V\T E^{-\e,*}_V@>>>E^{-\e,*}_V$ are the 
projections and the restriction of $(,)$ to $E^{\e,\l}_V\T E^{-\e,*}_V$ is denoted again by $(,)$. Thus it is 
enough to show that $\d'_{1!}\cl^{(,)}=\a_\l^*h'_!\cl^h$. We can assume that $V=\vv$. Then $E^{\e,\l}_V$ can be 
identified with $h'{}\i(\l)$, $E^{-\e,*}_V$ can be identified with $(\kk^*)^m$, $(,)$ can be identified with the 
map 
$$h'':h'{}\i(\l)\T(\kk^*)^m@>>>\kk,\qua ((x_1,\do,x_m),(y_1,\do,y_m))\m x_1y_1+\do+x_my_m,$$
$\d'_1$ can be identified with $\d'':h'{}\i(\l)\T(\kk^*)^m@>>>(\kk^*)^m$ (second projection) and $\a_\l$ can be 
identified with 
$$\a':(\kk^*)^m@>>>\kk^*,\qua (z_1,\do,z_m)\m\l z_1\do z_m.$$ 
Define $j:h'{}\i(\l)\T(\kk^*)^m@>>>(\kk^*)^m$ by 
$$j((x_1,\do,x_m),(y_1,\do,y_m))=(x_1y_1,\do,x_my_m).$$ 
It is enough to show that $\a'{}^*h'_! h^*=\d''_!h''{}^*$. 
Since $h''=hj$ it is enough to show that $\a'{}^*h'_! h^*=\d''_!j^*h^*$, or that $\a'{}^*h'_!=\d''_!j^*$. This 
follows from the fact that the diagram consisting of $\a',h',\d'',j$ is cartesian.

\subhead 2.3\endsubhead
Let $V\in\cc$. Let $\s\in\cp^{ap},t\in\NN$ be such that $|V|=|\s|+t\bold1$. Let $E^{-\e}_{V,t,\s}$ be the set of 
all $T\in E^{-\e}_V$ such that $\sha(\cs_T)=t$, $V_{T,\l}\cong\vv$ for $\l\in\cs_T$ and 
${}^0T\in E^{-\e}_{V_{T,0},\bsi}$ (${}^0T$ as in 1.1). Note that $E^{-\e}_{V,t,\s}$ is a locally closed, smooth,
irreducible, $G_V$-stable subvariety of $E^{-\e}_V$. We define a finite principal covering 
$\z:\tE^{-\e}_{V,t,\s}@>>>E^{-\e}_{V,t,\s}$ as follows. By definition, $\tE^{-\e}_{V,t,\s}$ is the set of all 
pairs $(T,\o)$ where $T\in E^{-\e}_{V,t,\s}$ and $\o$ is a total order on $\cs(T)$. The group of this covering is
$\fS_t$ (it acts freely in an obvious way on $\tE^{-\e}_{V,t,\s}$. The map $\z$ is $(T,\o)\m T$.

\subhead 2.4\endsubhead
Let $\hcp=\cp^{ap}\T\car$. For $\nu\in\NN^I$ let $\hcp_\nu=\{(\s,\r)\in\hcp;|\s|+\un\r\bold1=\nu\}$. We define a
map $\hcp_\nu@>>>\cp_\nu$ by $(\s,\r)\m\ti\s$ where $\ti\s(\b)=\s(\b)+\r(b-a)$ for any $\b=\ov{a,b}\in\ca_m$. 
Note that:

(a) {\it this map is a bijection.}
\nl
The inverse map is $\ti\s\m(\s,\r)$ where $\r\in\car$ is defined by 

$\r(n)=\min_{\b=\ov{a,b}\in\ca_m;b-a=n}\ti\s(\b)$
\nl
and $\s\in\cp^{ap}$ is defined by $\s(\b)=\ti\s(\b)-\r(b-a)$ for any $\b=\ov{a,b}\in\ca_m$. 

Let $Z$ be the set of all collections $\p=(\p^\l)_{\l\in\kk}$ where $\p^\l\in\car$ for $\l\in\kk^*$, $\p^0\in\hcp$
are such that $\p^\l=0$ for all but finitely many $\l\in\kk^*$. For $\nu\in\NN^I$ let 
$Z_\nu=\{\p\in Z;\sum_{\l\in\kk^*}\un{\p^\l}\bold1+|\p^0|=\nu\}$. For $\p\in Z$ and $\l\in\kk$ we define 
$\p_\l\in\car$ by $\p_\l=\p^\l$ if $\l\in\kk^*$, $\p_0=\r$ if $\l=0$ and $\p^0=(\s,\r)$. 

Assume that we are given a collection of bijections

(b) $\Ph_{n,\l}:\car_n@>\si>>\car_n$ ($\l\in\kk^*,n\in\NN$), $\Ps_{\nu'}:\cp_{\nu'}@>\si>>\hcp_{\nu'}$ 
($\nu'\in\NN^I$).
\nl
Such a collection exists by (a). We define a bijection $\Ph:\tZ@>\si>>Z$ by $(\th_\l)\m(\p_\l)$ where 
$\p_\l=\Ph_{n,\l}(\th_\l)$ for $\l\in\kk^*$ (here $n=\un{\th_\l}$) and $\p_0=\Ps_{\nu'}(\th_0)$ (here
$\nu'=|\th_0|$). This restricts for any $\nu\in\NN^I$ to a bijection 

(c) $\Ph_{(\nu)}:\tZ_\nu@>\si>>Z_\nu$.

\subhead 2.5\endsubhead
Let $V\in\cc$ and let $\nu=|V|$. Let $\p\in Z_\nu$. Let $z=\sum_{\l\in\kk^*}\un{\p_\l}$. We have
$\p^0=(\s,\p_0)\in\hcp$ where $\s\in\cp^{ap}$ and $|V|=|\s|+z\bold1$. We associate to $\p$ a finite unramified 
covering $\x:{}'E^{-\e}_{V,z,\s}@>>>E^{-\e}_{V,z,\s}$ as follows. By definition, ${}'E^{-\e}_{V,z,\s}$ is the set
of all pairs $(T,g)$ where $T\in E^{-\e}_{V,z,\s}$ and $g:\cs(T)@>>>\kk$ is a map such that 
$\sha(g\i(\l'))=\un{\p_{\l'}}$ for $l'\in\kk$. We define $\x$ by $(T,g)\m T$. Note that for 
$T\in E^{-\e}_{V,z,\s}$ we have 
$$\sha(\x\i(T))=z!(\prod_{\l'\in\kk}\un{\p_{\l'}}!)\i.$$ 
Let $\LL$ be the local system of rank $m^{z-\un{\p_0}}$ on ${}'E^{-\e}_{V,z,\s}$ whose stalk at $(T,g)$ is 
$\bxt_{\l\in\cs_T,g(\l)\ne0}\fK^m_{\l g(\l)}$. We define a finite principal covering 
$\z:{}''E^{-\e}_{V,z,\s}@>>>{}'E^{-\e}_{V,z,\s}$ as follows. By definition,
${}''E^{-\e}_{V,z,\s}$ is the set of all triples $(T,g,\o)$ where $(T,g)\in{}'E^{-\e}_{V,z,\s}$ and $\o$ is a 
collection of total orders on each of the sets $g\i(\l')$ ($\l'\in\kk$). The group of this covering is
$\cg:=\prod_{\l'\in\kk}\fS_{\un{\p_{\l'}}}$. (It acts freely in an obvious way on ${}''E^{-\e}_{V,z,\s}$.) Let 
$\tce$ be the local system on ${}'E^{-\e}_{V,z,\s}$ associated to the principal covering $\z$ and the irreducible
representation $\bxt_{\l'\in\kk}\Pi_{\p_{\l'}}$ of $\cg$. (The rank of $\tce$ is 
$\prod_{\l'\in\kk}N_{\p_{\l'}}$.) Let $\ce_\p=\x_!(\LL\ot\tce)$, a local system on $E^{-\e}_{V,z,\s}$ of rank 
$$z!(\prod_{\l'\in\kk}\un{\p_{\l'}}!)\i m^{z-\un{\p_0}}\prod_{\l'\in\kk}N_{\p_{\l'}}.$$
Hence the complex
$$IC(\ov{E^{-\e}_{V,z,\s}},\ce_\p)_{E^{-\e}_V}\tag a$$
is defined.

In the following theorem we assume that the characteristic of $\kk$ is sufficiently large. (The reason is 
explained in 2.9.) It is likely that this assumption is unnecessary.

\proclaim{Theorem 2.6} Let $V\in\cc$. Let $\nu=|V|$. 

(a) For any $\p\in Z_\nu$, the local system $\ce_\p$ on $E^{-\e}_{V,t,\s}$ (notation of 2.5) is irreducible.

(b) There exists a bijection $\tZ_\nu@>\si>>Z_\nu$ as in 2.4(c) such that if $\th\in\tZ_\nu$ and $\p\in Z_\nu$ is
the corresponding element under this bijection, then
$$\fF(IC(\ov{E^\e_{V,\th}},\bbq)_{E^\e_V}[?])\cong IC(\ov{E^{-\e}_{V,t,\s}},\ce_\p)_{E^{-\e}_V}[?].$$
\endproclaim
The proof will occupy much of the remainder of this section. We see that the complexes 2.5(a) for various 
$\p\in Z_\nu$ are (up to shift) exactly the antiorbital complexes on $E^{-\e}_V$.

\subhead 2.7\endsubhead
We preserve the setup of 1.5. Note that for $j\in[1,s]$, $\fF(\ck_j)$ is a $G_{U^j}$-equivariant perverse sheaf 
(up to shift) on $E^{-\e}_{U^j}$. Hence $\Ind(\bxt_j(\fF(\ck_j))\in\cd(E^{-\e}_V)$ is defined. We have:
$$\fF(\Ind(\bxt_j\ck_j))\cong\Ind(\bxt_j(\fF(\ck_j))[?]\text{ in }\cd(E^{-\e}_V).\tag a$$
This is proved in \cite{\LI, 5.4} assuming that $m\ge2$ and $s=2$; the proof for the case when $m\ge2$, $s\ne2$ is
similar to that for $m\ge2,s=2$ or can be deduced from the case $s=2$ by repetition. The proof for $m=1$ is
identical to that for $m\ge2$. (Alternatively, (a) can be deduced from Theorem A.2 in the Appendix.)

\subhead 2.8\endsubhead
We preserve the setup of 1.5 but we replace $\e$ by $-\e$. We fix $\l\in\kk^*$. Assume that $U^j\cong\vv$ for 
$j\in[1,s]$. Let $\a_{\l,j}:E^{-\e,*}_{U^j}@>>>\kk^*$ be the map $\a_\l$ of 1.1 with $V$ replaced by $U^j$. Let 
$\ck_j=IC(E^{-\e}_{U^j},\a_{\l,j}^*\fK^m)$. Let $S''$ be the set of all $(T,f)\in E''$ such that the
following holds: $f=(V^j)$ is such that for $j\in[1,s]$, $V^j/V^{j-1}\cong U^j$ and the element of
$E^\e_{V^j/V^{j-1}}$ induced by $T$ is in $E^{-\e,*}_{V^j/V^{j-1}}$. This is an open dense smooth irreducible 
subset of $E''$. Define $\a:S''@>>>(\kk^*)^s$ by $\a(T,f)=(a_1\l,a_2\l,\do,a_s\l)$ where $T^m=a_j$ on 
$V^j/V^{j-1}$. Let $\cl''=\a^*((\fK^m)^{\bxt s})$, a local system on $S''$. We have $\ck''=IC(E'',\cl'')$. Let $S$
be the set of all $T\in E^{-\e}_V$ such that $\sha(\cs_T)=s$ and $V_{T,\l}\cong\vv$ for $\l\in\cs_T$. (Thus, 
$S=E^{-\e}_{V,s,0}$, see 2.3.) This is an open dense subset of $E^{-\e}_V$. Also $p_3\i(S)$ is an open dense 
(hence irreducible) subset of $S''$; let $p'_3:p_3\i(S)@>>>S$ be the restriction of $p_3$. We have an isomorphism 
$u:\tE^{-\e}_{V,s,0}@>\si>>p_3\i(S)$ (notation of 2.3) given by $(T,\o)\m(T,f)$ where $f=(V^j)$ is defined by 
$V^1=V_{T,\l_1}$, $V^2=V_{T,\l_1}\op V_{T,\l_2}$, etc. Here $\l_1,\l_2,\do,\l_s$ are the elements of $\cs(T)$ 
arranged in the order given by $\o$. Under the isomorphism $u$, the map $p'_3$ becomes the map $\z$ in 2.3 hence 
is a finite principal covering whose group is $\fS_s$. Thus $p'_{3!}\bbq\cong\op_{\r\in\car_s}(\cl^\r)^{\op N_\r}$
as local systems on $S$ where $\cl^\r$ is the local system
on $S$ associated to $p'_3$ and the irreducible representation $\Pi_\r$ of $\fS_s$. Note that the local systems 
$\cl^\r,\r\in\car_s$, are irreducible and mutually nonisomorphic (by the irreducibility of $p_3\i(S)$). Define 
$\bar\a:S@>>>(\kk^*)^s/\fS_s$ by $\bar\a(T)=(a_1\l,a_2\l,\do,a_s\l)$ (unordered) where $a_j\in\kk^*$ are the 
scalars such that $V_{T,a_j}\ne0$. Let $\bar\LL$ be the local system on $(\kk^*)^s/\fS_s$ whose inverse image
under the obvious map $(\kk^*)^s@>>>(\kk^*)^s/\fS_s$ is $(\fK^m)^{\bxt s}$. Let $\LL=\bar\a^*(\bar\LL)$, a local
system on $S$. We have $\cl''|_{p_3\i(S)}=p'_3{}^*\LL$. Hence
$$\Ind(\bxt_j\ck_j)|_S=p'_{3!}(p'_3{}^*\LL)=\LL\ot p'_{3!}\bbq=\op_{\r\in\car_s}\LL\ot(\cl^\r)^{\op N_\r}.$$
By 2.2, for $j\in[1,s]$, our $\ck_j$ is (up to shift) the D-F transform of $\ck_j$ in 1.6. Using 2.7 we 
deduce that our $\Ind(\bxt_j\ck_j)$ is (up to shift) the D-F transform of $\Ind(\bxt_j\ck_j)$ of 1.6. Hence
$$\Ind(\bxt_j\ck_j)\cong\op_{\r\in\car_s}\fF(P''_\r)^{\op N_\r}[d]$$
with $P''_\r$ as in 1.6 and $d\in\ZZ$. Note that $P''_\r$ are simple mutually nonisomorphic perverse sheaves on
$E^{-\e}_V$. Restricting to $S$ we obtain
$$\op_{\r\in\car_s}\LL\ot(\cl^\r)^{\op N_\r}\cong\op_{\r\in\car_s}\fF(P''_\r)|_S^{\op N_\r}[d]\tag a$$
in $\cd(S)$. It follows that for each $\r$, $\fF(P''_\r)|_S[d]$ is a local system on $S$ (possibly zero). Let 
$\car'_s$ be the set of all $\r\in\car_s$ such that $\fF(P''_\r)|_S[d]\ne0$. Since $S$ is open dense in 
$E^{-\e}_V$, $\fF(P''_\r)|_S[d]$ ($\r\in\car'_s$) are irreducible mutually nonisomorphic local system on $S$. Thus
in (a), the right hand side is a direct sum of $\sum_{\r\in\car'_s}N_\r$ irreducible local systems while the left
hand side is a direct sum of $\sum_{\r\in\car_s}N_\r$ nonzero local systems. This forces $\car_s=\car'_s$ and 
that each $\LL\ot\cl^\r$ is irreducible. If $\r\ne\r'$ then $\LL\ot\cl^\r\not\cong\LL\ot\cl^{\r'}$; otherwise the
number of nonisomorphic irreducible local systems which appear in the left hand side of (a) would be
$<\sha(\car_s)$ while the analogous number for the right hand side of (a) would be equal to $\sha(\car_s)$. We see
that there is a unique permutation $\r\m\r'$ of $\car_s$ such that $\fF(P''_\r)|_S[d]\cong\LL\ot\cl^{\r'}$ for any
$\r\in\car_s$. It follows that $\fF(P''_\r)\cong IC(E^{-\e}_V,\LL\ot\cl^{\r'})_{E^{-\e}_V}[?]$. Thus the D-F 
transforms of orbital complexes on $E^\e_V$ with support contained in $E^{\e,\l}_V$ are exactly the complexes of 
the form $IC(E^{-\e}_V,\LL\ot\cl^{\r'})_{E^{-\e}_V}[?]$ for various $\r\in\car_s$. The bijection 
$\car_s@>>>\car_s$, $\r\m\r'$ is denoted by $\Ph_{s,\l}$.

\subhead 2.9\endsubhead
Let $V\in\cc$. Combining two results \cite{\LI, 10.14}, \cite{\LII, 5.9} in the theory of canonical bases, we see
that 

(a) {\it the collection of simple perverse sheaves 
$IC(\ov{E^\e_{V,\bsi}},\bbq)_{E^\e_V}[?]$, ($\s\in\cp^{ap}_{|V|}$)
is mapped bijectively by $\fF$ onto the collection of simple perverse sheaves 

$IC(\ov{E^{-\e}_{V,\bsi}},\bbq)_{E^{-\e}_V}[?]$, ($\s\in\cp^{ap}_{|V|}$).}
\nl
(Both collections index the canonical basis in degree $|V|$ of a certain algebra associated to affine $SL_2$. Note
that in \cite{\LII} the arguments are in characteristic $0$ and they imply what we need only in sufficiently large
characteristic.)

It follows that there exists a bijection $\s\m\s^*$, $\cp^{ap}_{|V|}@>\si>>\cp^{ap}_{|V|}$ such that
$$\fF(IC(\ov{E^\e_{V,\bsi}},\bbq)_{E^\e_V}[?])\cong IC(\ov{E^{-\e}_{V,\bsi^*}},\bbq)_{E^{-\e}_V}[?]\tag b$$
for any $\s\in\cp^{ap}_{|V|}$.

\subhead 2.10\endsubhead
Let $V\in\cc$. Let $\ck=IC(\{0\},\bbq)_{E^\e_V}$, $\ck'=\bbq=IC(E^{-\e}_V,\bbq)$. From the definitions we have 
$\fF(\ck)=\ck'[n]$ where $n=\dim E^{-\e}_V$. 

\subhead 2.11\endsubhead
We preserve the setup of 1.5 but we replace $\e$ by $-\e$. Assume that $s=t+1$ where $t\in\NN$. Assume that
$U^j,\s$ are as in 1.7. Let $\s^*\in\cp^{ap}_{|U^s|}$ be as in 2.9. For $j\in[1,t]$ let 
$\ck_j=\bbq\in\cd(E^{-\e}_{U^j})$. Let $\ck_s=IC(\ov{E^{-\e}_{U^s,\bsi^*}},\bbq)_{E^{-\e}_{U^s}}$, ($\bsi^*$ as in
1.3). Let $S''$ be the set of all $(T,f)\in E''$ such that the following holds: $f=(V^j)$ is such that 
$V^j/V^{j-1}\cong U^j$ for $j\in[1,s]$ and the element of $E^{-\e}_{V/V^{s-1}}$ induced by $T$ is in 
$E^{-\e}_{V/V^{s-1},\bsi^*}$. This is a locally closed smooth irreducible subset of $E''$. In our case the 
complex $\ck''$ of 1.5 is $\ck''=IC(\bS'',\bbq)_{E''}$. Let $S=E^{-\e}_{V,t,\s^*}$, see 2.3. Now $p_3\i(S)$ is an 
open dense (hence irreducible) subset of $S''$; let $p'_3:p_3\i(S)@>>>S$ be the restriction of $p_3$. We have an 
isomorphism $u:\tE^{-\e}_{V,t,\s^*}@>\si>>p_3\i(S)$ (notation of 2.3) given by $(T,\o)\m(T,f)$ where $f=(V^j)$ is 
defined by $V^1=V_{T,\l_1}$, $V^2=V_{T,\l_1}\op V_{T,\l_2}$, etc. Here $\l_1,\l_2,\do,\l_t$ are the elements of 
$\cs(T)$ arranged in the order given by $\o$. Under the isomorphism $u$, the map $p'_3$ becomes the map $\z$ in 
2.3 hence is a finite principal covering whose group is $\fS_t$. Thus 
$p'_{3!}\bbq\cong\op_{\r\in\car_t}(\cl^\r)^{\op N_\r}$ as local systems on $S$ where $\cl^\r$ is the local system
on $S$ associated to $p'_3$ and the irreducible representation $\Pi_\r$ of $\fS_t$. Note that the local systems 
$\cl^\r,\r\in\car_t$, are irreducible and mutually nonisomorphic (by the irreducibility of 
$p_3\i(S)$). Since $p_3$ is proper, we have $p_3(\bS'')=p_3(\ov{p_3\i(S)})=\ov{p_3(p_3\i(S))}=\bS$ where $\bar{}$
denotes closure. Hence $\supp(p_{3!}\ck'')\sub\bS$. 

By 2.9(b), 2.10, for $j\in[1,s]$, our $\ck_j$ is (up to shift) the D-F transform of $\ck_j$ in 1.7. Using 2.7 we 
deduce that our $\Ind(\bxt_j\ck_j)$ is (up to shift) the D-F transform of $\Ind(\bxt_j\ck_j)$ of 1.7. Hence 
$\Ind(\bxt_j\ck_j)=p_{3!}\ck''$ is isomorphic to $\op_{e\in[1,M]}\fF(P_e)[d'_e]$ where $P_e$ are as in 1.7 and 
$d'_e\in\ZZ$. Note that $\fF(P_e)$ are simple perverse sheaves on $E^{-\e}_V$ with support contained in $\bS$.
Restricting to $S$ we obtain
$$\op_{\r\in\car_t}(\cl^\r)^{\op N_\r}\cong\op_{e\in[1,M]}\fF(P_e)|_S[d'_e]$$
in $\cd(S)$. It follows that for any $e$, $\fF(P_e)|_S[d'_e]$ is a local system on $S$ (possibly zero). Since 
$\supp(\fF(P_e))\sub\bS$ we see that $\fF(P_e)|_S[d'_e]$ is either an irreducible local system on $S$ or is $0$. 
We deduce that $\fF(P_e)|_S$ is either $0$ or is isomorphic to $\cl^\r[-d'_e]$ for a well defined $\r\in\car_t$. 
Hence either $\supp(\fF(P_e))\sub\bS-S$ or $\fF(P_e)\cong IC(\bS,\cl^\r)_{E^{-e}_V}[?]$. Conversely, we see that 
for any $\r\in\car_t$ we have $IC(\bS,\cl^\r)_{E^{-e}_V}[?]\cong\fF(P_e)$ for some $e$. Let $\cx_{t,\s^*}$ be the 
set of isomorphism classes of simple perverse sheaves on $E^{-\e}_V$ that are isomorphic to $\fF(P_e)$ for some 
$e\in[1,M]$. Let $\cx'_{t,\s^*}$ be the set of isomorphism classes of simple perverse sheaves on $E^{-\e}_V$ that 
are of the form $IC(\bS,\cl^\r)_{E^{-e}_V}[?]$ for some $\r\in\car_t$. We see that

(a) $\sha(\car_t)=\sha(\cx'_{t,\s^*})\le\sha(\cx_{t,\s^*})$.
\nl
Let $\cy'_{t,\s}$ be the subset of $\cy_{t,\s}$ (see 1.7) which corresponds to the subset $\cx'_{t,\s^*}$ under 
the bijection 

(b) $\cy_{t,\s}@>\si>>\cx_{t,\s^*}$.
\nl
induced by $\fF$. Then $\fF$ defines a bijection

(c) $\cy'_{t,\s}@>\si>>\cx'_{t,\s^*}$.
\nl
We show:

(d) if $t,\tit\in\NN$ and $\s,\ti\s\in\cp^{ap}$ are such that $|V|=|\s|+t\bold1=|\ti\s|+\tit\bold1$ and 
$(t,\s)\ne(\tit,\ti\s)$, then $\cx'_{t,\s^*}\cap\cx'_{\tit,\ti\s^*}=\em$. 
\nl
It is enough to show that $E^{-\e}_{V,t,\s^*}\cap E^{-\e}_{V,\tit,\ti\s^*}=\em$. This is clear from the 
definition.

Applying the inverse of $\fF$ we obtain:

(e) in the setup of (d), we have $\cy'_{t,\s}\cap\cy'_{\tit,\ti\s}=\em$. 
\nl
Let $\nu=|V|$. Let $X$ be the set of all $(t,\ti\s)\in\NN\T\cp^{ap}$ such that $|\ti\s|+t\bold1=\nu$. We have
$$\align&\sha(\cp_\nu)\ge\sha(\cup_{(t,\ti\s)\in X}\cy_{t,\ti\s})\ge\sha(\cup_{(t,\ti\s)\in X}\cy'_{t,\ti\s})
=\sum_{(t,\ti\s)\in X}\sha(\cy'_{t,\ti\s})\\&=\sum_{(t,\ti\s)\in X}\sha(\cx'_{t,\ti\s^*})=
\sum_{(t,\ti\s)\in X}\sha(\car_t)=\sha(\hcp_\nu)=\sha(\cp_\nu).\tag f\endalign$$
The first $\ge$ follows from the fact that each $\cy_{t,\s}$ consists of simple perverse sheaves of the form
$IC(\ov{E^\e_{V,\bsi}},\bbq)_{E^\e_V}[?]$ with $\s\in\cp_\nu$. The second $\ge$ follows from the inclusion 
$\cy'_{t,\ti\s}\sub\cy_{t,\ti\s}$. The first $=$ follows from (e). The second $=$ follows from (c). The third $=$ 
follows from (b). The fourth $=$ follows from definitions. The fifth $=$ follows from 2.4(a). It follows that 
each inequality in (f) is an equality. 

In particular we see that any simple perverse sheaf of the form $IC(\ov{E^\e_{V,\ov{\s_1}}},\bbq)_{E^\e_V}[?]$ 
with $\s_1\in\cp_\nu$ belongs to $\cy'_{t,\s}$ for a unique $(t,\s)\in X$. Hence any simple perverse sheaf of the
form $\fF(IC(\ov{E^\e_{V,\ov{\s_1}}},\bbq)_{E^\e_V}[?])$ with $\s_1\in\cp_\nu$ belongs to $\cx'_{t,\s^*}$ for a
unique $(t,\s)\in X$. In particular it is of the form $IC(\ov{E^{-e}_{V,t,\s^*}},\cl^\r)_{E^{-e}_V}[?]$ for a 
unique $(\s^*,\r)\in\hcp_{|V|}$, $t=\un\r$.

Since $\sha(\cp_\nu)=\sha(\hcp_\nu)$, see 2.4, we see that there exists a unique bijection 
$\Ps_\nu:\cp_\nu@>\si>>\hcp_\nu$, $\s'\m\hat\s'$, such that
$$\fF(IC(\ov{E^\e_{V,\bsi'}},\bbq)_{E^\e_V}[?])\cong IC(\ov{E^{-e}_{V,t,\s^*}},\cl^\r)_{E^{-e}_V}[?]\tag g$$
for any $\s'\in\cp_\nu$; here $\hat\s'=(\s^*,\r)$, $t=\un\r$.

\subhead 2.12\endsubhead
We preserve the setup of 1.5 but we replace $\e$ by $-\e$. Assume that $s=t+1$, $t\in\NN$. Let 
$\l_1,\l_2,\do,\l_s$, $\th$, $U^j$ be as in 1.8. Let $\nu=|V|$. Let $\Ph_{(\nu)}:\tZ_\nu@>\si>>Z_\nu$ be the 
bijection associated in 2.4 to the collection of bijections $\Ph_{n,\l}$ (as in 2.8) and $\Ps_{\nu'}$ (as in 
2.11). Let $\p=\Ph_{(\nu)}(\th)\in Z_\nu$. For $j\in[1,t]$ let $d_j=\un{\p_{\l_j}}$. Let 
$S_j=E^{-\e}_{U^j,d_j,0}$, see 2.3. Let $\LL_j$ be the local system on $S_j$ defined as $\LL$ in 2.8 (with 
$V,s,\l$ replaced by $U^j,d_j,\l_j$). Let $\cl^{\p_{\l_j}}$ be the local system on $S_j$ defined as $\cl^{\r'}$ in
2.8 (with $V,s,\l,\r$ replaced by $U^j,d_j,\l_j,\th_{\l_j}$). Let 
$\ck_j=IC(E^{-\e}_{U^j},\LL_j\ot\cl^{\p_{\l_j}})_{E^{-\e}_{U^j}}$. By 2.8, this $\ck_j$ is (up to shift) the D-F 
transform of the $\ck_j$ in 1.8. We write $\p^{\l_s}=\p^0$ as $(\s,\r)\in\cp^{ap}\T\car$. Let $d_s=\un\r$. Let 
$S_s=E^{-\e}_{U^s,d_s,\s}$. Let $\cl^{\p_0}$ be the local system on $S_s$ defined as $\cl^\r$ in 2.11 (with 
$V,t,\r$ replaced by $U^s,d_s,\p_0$). Let $\ck_s=IC(\ov{E^{-\e}_{U^s,d_s,\s}},\cl^{\p_0})_{E^{-\e}_{U^s}}$. By 
2.11, this $\ck_s$ is (up to shift) the D-F transform of the $\ck_s$ in 1.8. 

Let $S''$ be the set of all $(T,f)\in E''$ such that the following holds: $f=(V^j)$ is such that 
$V^j/V^{j-1}\cong U^j$ for $j\in[1,s]$ and $T$ is such that the element $T_j$ of $E^{-\e}_{V^j/V^{j-1}}$ induced 
by $T$ is in $E^{-\e}_{V^j/V^{j-1},d_j,0}$ (if $j\in[1,t]$) and in $E^{-\e}_{V/V^s,d_s,\s}$ (if $j=s$). This is a
locally closed smooth irreducible subset of $E''$. We define a finite principal covering $\z'':\tS''@>>>S''$ as 
follows. By definition, $\tS''$ is the set of all triples $(T,f,\o)$ where $(T,f)\in S''$ and $\o=(\o_j)$ is a 
collection of total orders $\o_j$ on $\cs(T_j)$ ($j=1,\do,s$). The group of this covering is 
$\cg=\prod_{j\in[1,s]}\fS_{d_j}$ (it acts freely in an obvious way on $\tE^{-\e}_{V,t,\s}$). The map $\z''$ is 
$(T,f,\o)\m(T,f)$. Let $\tce$ be the local system on $S''$ associated to the principal covering $\z''$ and the 
irreducible representation $\bxt_{\l'\in\kk}\Pi_{\p_{\l'}}$ of $\cg$. Let $\LL''$ be the local system on $S''$ 
whose stalk at $(T,f)$ is $\bxt_{j\in[1,t]}\bxt_{\l\in\cs_{T_j}}\fK^m_{\l\l_j}$. In our case the complex $\ck''$ 
of 1.5 is $\ck''=IC(\bS'',\LL''\ot\tce)_{E^{-e}_V}$. From 2.7 we see that our $\Ind(\bxt_j\ck_j)=p_{3!}\ck''$ is 
(up to shift) the D-F transform of $\Ind(\bxt_j\ck_j)=IC(\ov{E^\e_{V,\th}},\bbq)_{E^\e_V}$ of 1.8.

Let $z=\sum_{j\in[1,s]}\un{\p_{\l_j}}$. Let $S=E^{-\e}_{V,z,\s}$, see 2.3. Now $p_3\i(S)$ is an open dense subset
of $S''$; let $p'_3:p_3\i(S)@>>>S$ be the restriction of $p_3$. We have an isomorphism 
$u:{}'E^{-\e}_{V,z,\s}@>\si>>p_3\i(S)$ (see 2.5) given by $(T,g)\m(T,f)$ where $f=(V^j)$ is defined by 
$$V^1=\op_{\l\in g\i(\l_1)}V_{T,\l},\qua V^2=\op_{\l\in g\i(\l_1)\cup g\i(\l_2)}V_{T,\l},\text{ etc}.$$
Under the isomorphism $u$, the map $p'_3$ becomes the map $\x$ of 2.5. Since $p_3$ is proper, we have 
$p_3(\bS'')=p_3(\ov{p_3\i(S)})=\ov{p_3(p_3\i(S))}=\bS$ where $\bar{}$ denotes closure. Hence 
$\supp(p_{3!}\ck'')\sub\bS$ so that $p_{3!}\ck''$ is completely determined by its restriction to $S$. Since 
$p_{3!}\ck''[?]$ is the D-F transform of a simple perverse sheaf, it is itself a simple perverse sheaf. We have 
$$(p_{3!}\ck'')|_S=p'_{3!}((\LL''\ot\tce)_{p_3\i(S)})=\x_!(\LL\ot\tce)=\ce_\p$$
(notation of 2.5). We see that $p_{3!}\ck''=IC(\ov{E^{-\e}_{V,z,\s}},\ce_\p)_{E^{-e}_V}$. Thus 
$IC(\ov{E^{-\e}_{V,z,\s}},\ce_\p)_{E^{-e}_V}$ is (up to shift) the D-F transform of 
$IC(\ov{E^\e_{V,\th}},\bbq)_{E^\e_V}$. This forces $\ce_\p$ to be an irreducible local system, proving 2.6(a); 
2.6(b) follows as well.

\subhead 2.13\endsubhead
We state a converse of 2.9(a). Let $V\in\cc$. Let $\nu=|V|$. Let $\ck$ be a $G_V$-equivariant simple perverse 
sheaf on $E^\e_V$. 

(a) {\it The following three conditions on $\ck$ are equivalent:

(i) $\supp(\ck)$ is contained in the nilpotent variety and $\supp(\fF(\ck))$ is contained in the nilpotent 
variety;

(ii) $\ck$ is biorbital.

(iii) $\ck\cong IC(\ov{E^\e_{V,\bsi}},\bbq)_{E^\e_V}[?]$ for some $\s\in\cp^{ap}_\nu$.}
\nl
(This gives a new characterization of the perverse sheaves which constitute the canonical basis \cite{\LI} 
associated to a cyclic quiver.)

Let $S_1$ (resp. $S_2$ or $S_3$) be the set of $\ck$ (up to isomorphism) as in (i) (resp. as in (ii) or (iii)). 
From 2.9(a) we see that $S_3\sub S_1$. From 2.11 we see that 

$\sha(S_1)=\sha((\s^*,\r)\in\hcp_\nu;\un\r=0\}=\sha(\cp^{ap}_\nu)=\sha(S_3)$.
\nl
It follows that $S_1=S_3$.

Clearly, $S_1\sub S_2$. Assume now that $\ck\in S_2$. If $\supp(\ck)$ is the closure of a non-nilpotent
orbit then from 2.12 we see that the support of $\fF(\ck)$ is a closure of a subvariety of the form
$E^{-\e}_{V,z,\s}$ where $z>0$, $\s\in\cp^{ap}$; in particular, the support of $\fF(\ck)$ is not the closure of a
single orbit so that $\fF(\ck)$ is not orbital, a contradiction. Thus, $\supp(\ck)$ is the closure of a nilpotent
orbit. The same argument shows that if $\supp(\fF(\ck))$ is the closure of a non-nilpotent orbit then 
$\fF(\fF(\ck))$ is not orbital hence $\ck$ is not orbital, a contradiction. Thus, $\supp(\fF(\ck))$ is the closure
of a nilpotent orbit. We see that $\ck\in S_1$. Thus $S_2\sub S_1$ hence $S_1=S_2$. This proves (a).

\subhead 2.14\endsubhead
Let $V\in\cc$. Let $\nu=|V|$. Assume that we are given an $\FF_q$-structure on each $V_i$. Then $E^\e_V$ inherits
an $\FF_q$-structure and for each $\s\in\cp_\nu$, the subset $E^\e_{V,\bsi}$ of $E^\e_V$ is defined over $\FF_q$.
Let $U=E^\e_V(\FF_q)$. Then the dual space $U^*$ can be identified with $E^{-\e}_V(\FF_q)$. For $\s\in\cp_\nu$ let
$U_\s=E^\e_{V,\bsi}(\FF_q)$, $U^*_\s=E^{-\e}_{V,\bsi}(\FF_q)$; we define $f_\s:U@>>>\bbq$, $f'_\s:U^*@>>>\bbq$ as
follows. If $x\in U\cap\ov{E^\e_{V,\bsi}}$ (resp. $x'\in U^*\cap\ov{E^{-\e}_{V,\bsi}}$) we define $f_\s(x)$ (resp.
$f'_\s(x')$) as the alternating sum of the traces of the Frobenius map on the stalks of the cohomology sheaves of
$IC(\ov{E^\e_{V,\bsi}},\bbq)$ (resp. $IC(\ov{E^{-\e}_{V,\bsi}},\bbq)$) at $x$ (resp. $x'$); if 
$x\in U-\ov{E^\e_{V,\bsi}}$ (resp. $x'\in U^*-\ov{E^{-\e}_{V,\bsi}}$) we set $f_\s(x)=0$ (resp. $f'_\s(x')=0$). 
Note that $(f_\s)_{\s\in\cp_\nu}$ (resp. $(f'_\s)_{\s\in\cp_\nu}$) is a $\bbq$-basis of the vector space $\cf$ 
(resp. $\cf'$) of functions $U@>>>\bbq$ (resp. $U^*@>>>\bbq$) which are constant on the orbits of $G_V(\FF_q)$ and
vanish on non-nilpotent elements. 

For any $(\s,\r)\in\hcp_\nu$ we define a function $f''_{\s,\r}:U^*@>>>\bbq$ as follows. Let $t=\un\r$. If 
$x\in U^*\cap\ov{E^{-\e}_{V,t,\s}}$ we define $f''_{\s,\r}(x)$ as the alternating sum of the traces of the 
Frobenius map on the stalks of the cohomology sheaves of $IC(\ov{E^{-\e}_{V,t,\s}},\cl^\r)$ at $x$; if 
$x\in U^*-\ov{E^{-\e}_{V,t,\s}}$, we set $f''_{\s,\r}(x)=0$. We have a partition 
$E^{-\e}_{V,t,\s}(\FF_q)=\sqc_gE^{-\e}_{V,t,\s}(\FF_q)_g$ ($g$ runs through the conjugacy classes in $\fS_t$) 
where $E^{-\e}_{V,t,\s}(\FF_q)_g$ is the set of all $T\in E^{-\e}_{V,t,\s}(\FF_q)$ such that the action of 
Frobenius on $\cs_T$ is by a permutation of type $g$. Clearly each $E^{-\e}_{V,t,\s}(\FF_q)_g$ is nonempty. From 
the definitions for any $T\in E^{-\e}_{V,t,\s}(\FF_q)$ we have $f''_{\s,\r}(T)=\tr(g,\Pi_\r)$ where $g$ is defined
by $T\in E^{-\e}_{V,t,\s}(\FF_q)_g$.

From 2.11(g) we see that there exists a bijection $\cp_\nu@>\si>>\hcp_\nu$, $\s'\m\hat\s'$, such that for any 
$\s'\in\cp_\nu$ we have $\hf_{\s'}=b_{\s'}f''_{\s^*,\r}$ where $(\s^*,\r)=\hat\s'$ and $b_{\s'}\in\bbq^*$. 
Moreover from 2.9(b) we see that we have $\s'\in\cp^{ap}_\nu$ if and only if $(\s^*,\r)=\hat\s'$ satisfies 
$\un\r=0$.

Let $f\in\cf$ be such that $\hf\in\cf'$. We can write uniquely $f=\sum_{\s\in\cp_\nu}a_\s f_\s$ where 
$a_\s\in\bbq$. Assume that $a_\s\ne0$ for some $\s\in\cp_\nu-\cp^{ap}_\nu$. Applying Fourier transform we obtain 
$\hf=\sum_{(\s,\r)\in\hcp_\nu}a'_{\s,\r}f''_{\s,\r}$ where $a'_{\s,\r}\in\bbq$ and $a'_{\s_0,\r_0}\in\bbq^*$ for 
some $(\s_0,\r_0)\in\hcp_\nu$ with $\un\r_0>0$. We can assume in addition that 
$\dim E^{-\e}_{V,\un\r_0,\s_0}\ge\dim E^{-\e}_{V,\un\r,\s}$ for any $(\s,\r)\in\hcp_\nu$ such that 
$a'_{\s,\r}\ne0$, $\un\r>0$. Let $t=\un\r_0$. Then for any $(\s,\r)\in\hcp_\nu$ such that $a'_{\s,\r}\ne0$ we have
$\ov{E^{-\e}_{V,\un\r,\s}}\cap E^{-\e}_{V,\un\r_0,\s_0}=\em$ unless $\un\r=t$, $\s=\s_0$. Since $\hf\in\cf'$ we 
have $\hf|_{E^{-\e}_{V,\un\r_0,\s_0}}=0$. Hence 

$\sum_{(\s,\r)\in\hcp_\nu}a'_{\s,\r}f''_{\s,\r}|_{E^{-\e}_{V,t,\s_0}}=0$
\nl
and 

$\sum_{\r\in\car_t}a'_{\s_0,\r}f''_{\s_0,\r}|_{E^{-\e}_{V,t,\s_0}}=0$,

$\sum_{\r\in\car_t}a'_{\s_0,\r}\tr(g,\Pi_\r)=0$
\nl
for any conjugacy class $g$ in $\fS_t$. It follows that $a'_{\s_0,\r}=0$ for any $\r\in\car_t$, contradicting  
$a'_{\s_0,\r_0}=0$. This contradiction shows that $a_\s=0$ for any $s\in\cp_\nu-\cp^{ap}_\nu$. Thus $f$ is a 
linear combination of the functions $f_\s,(\s\in\cp^{ap}_\nu)$. Conversely, if $\s\in\cp^{ap}_\nu$ then from 
2.9(b) we see that $\hf_\s=cf'_{\s^*}$ where $\s^*\in\cp^{ap}_\nu$ and $c\in\bbq^*$. Thus we have $\hf_s\in\cf'$.
Thus any linear combination $f$ of the functions $f_\s,(\s\in\cp^{ap}_\nu)$ satisfies $f\in\cf$, $\hf\in\cf'$ and
we have the following result.

(a) {\it The functions $f_\s,(\s\in\cp^{ap}_\nu)$ form a basis of the vector space of all functions $f\in\cf$ such
that $\hf\in\cf'$.}

\subhead 2.15\endsubhead
The results in 2.11 suggest a way to organize the nilpotent $K$-orbits in $\fg_\z$ (in the context of 0.3) with
$\z\ne1$. (Here we assume for simplicity that all nilpotent elements in $\fg_\z$ have connected isotropy group in
$K$ but a similar picture should hold in general.) Namely, each nilpotent $K$-orbit should be attached to a 
$\vt$-stable Levi subgroup $L$ of $G$ up to $K$-conjugacy (such that $\vt$ acts on $L_{der}$, the derived group of
$L$, as an inner automorphism) and a biorbital complex $\ck$ on the $\z$-part of the Lie algebra of $L_{der}$.
Moreover, the nilpotent $K$-orbits corresponding to a given $(L,\ck)$ should be such that the D-F transforms of
the corresponding orbital complexes have the same support and they should be indexed by something similar to the 
irreducible representations of a Weyl group. Thus something like the "generalized Springer correspondence" should
hold even though the small and semismall maps in the usual theory are missing in general.

\head 3. Further examples of antiorbital complexes\endhead
\subhead 3.1\endsubhead
Let $G,\fg,\vt,K,\k$ be as in 0.3. We assume that the characteristic of $\kk$ is sufficiently large. For 
$\z\in\kk^*$ let $\fg_\z,\k_\z$ be as in 0.3. For any subspace $V$ of $\fg$ we set $V^\pe=\{x\in\fg;\k(x,V)=0\}$.

We say that $\vt:G@>>>G$ is {\it inner} if there exists a semisimple element $g_0\in G$ such that 
$\vt(g)=g_0gg_0\i$ for all $g\in G$. 

For any parabolic subgroup $P$ of $G$ we denote by $U_P$ the unipotent radical of $P$ and by $\uP,\un{U}_P$ the
Lie algebras of $P,U_P$. Assume now that $\vt(P)=P$. Then $\uP=\op_{\z\in\kk^*}\uP_\z$, 
$\un{U}_P=\op_{\z\in\kk^*}\un{U}_{P,\z}$ where $\uP_\z=\uP\cap\fg_\z$, $\un{U}_{P,\z}=\un{U}_P\cap\fg_\z$.

Let $\z\in\kk^*$. Assume that $\vt$ is inner. Let $X$ be a $K$-orbit in the variety of $\vt$-stable parabolic 
subgroups of $G$ (or equivalently a connected component of that variety). Note that $X$ is a projective variety. 
Let $\tX_\z=\{(x,P);P\in X,x\in\uP_\z\}$, $\tX'_\z=\{(x,P);P\in X,x\in\un{U}_{P,\z}\}$. Define
$\p_\z:\tX_\z@>>>\fg_\z$, $\p'_\z:\tX'_\z@>>>\fg_\z$ by $(x,P)\m x$. We show:
$$\fF(\p_{\z!}(\bbq))=\p'_{\z\i!}(\bbq))[?].\tag a$$
Here we view $\fF$ as a functor $\cd(\fg_\z)@>>>\cd(\fg_{\z\i})$. We have a commutative diagram
$$\CD
\tX_\z@<b<<      \Xi@>c>>\bXi  \\
@V\p_\z VV     @VrVV      @VdVV  \\
\fg_\z@<s<<\fg_\z\T\fg_{\z\i}@>t>>\fg_{\z\i}
\endCD$$
where 

$\Xi=\{(x,x',P);P\in X,x\in\uP_\z,x'\in\fg_{\z\i}\}$, $\bXi=\{(x',P);x'\in\fg_{\z\i},P\in X\}$, 
\nl
$s,t$ are the obvious projections and $b,c,r,d$ are the obvious maps. Let $\tit=tr:\Xi@>>>\fg_{\z\i}$. Let 
$\k':\Xi@>>>\kk$ be $(x,x',P)\m\k_\z(x,x')$. We have 
$$\fF(\p_{\z!}(\bbq))=t_!(r_!(\bbq)\ot\cl^{\k_\z})[?]=\tit_!(\cl^{\k'})[?]=d_!c_!(\cl^{\k'})[?].$$
We have a partition $\Xi=\Xi_0\cup\Xi_1$ where $\Xi_0$ is defined by the condition that $x'\in\un{U}_P$. Let 
$c_1:\Xi_1@>>>\bXi$ be the restriction of $c$. We show that 
$$c_{1!}(\cl^{\k'})=0.\tag b$$
The fibre of $c_1$ at $(x',P)$ is $\uP_\z$. The restriction of $\k'$ to this fibre is the linear map 
$x\m\k(x,x')$. It is enough to show that this linear map is not identically zero. (Assume that $\k(x,x')=0$ for 
any $x\in\uP_\z$ that is, $x'\in(\uP_\z)^\pe$. Since $x'\in\fg_{\z\i}$ we have automatically 
$x'\in(\uP_{\z'})^\pe$ for $\z'\ne\z$ hence $x'\in\uP^\pe$. Hence $x'\in\un{U}_P$ a contradiction.) This proves 
(b).

From (b) we deduce
$$d_!c_!(\cl^{\k'})=d_!c_!j_!j^*\cl^{\k'}=t_{0!}(\cl^{\k'}|_{\Xi_0})$$
where $j:\Xi_0@>>>\Xi$ is the inclusion and $t_0:\Xi_0@>>>\fg_\z$ is $(x,x',P)\m x'$. Now $t_0$ is a composition 
$\Xi_0@>>>\tX'_{\z\i}@>\p'_{\z\i}>>\fg_{\z\i}$ (the first map, $(x,x',P)\m(x',P)$, is an affine bundle whose fibre
at $(x',P)$ is isomorphic to $\uP_\z$, a vector space of constant dimension as $P$ runs through $X$ which is a 
$K$-orbit.) If $(x,x',P)\in\Xi_0$ then $\k'(x,x',P)=\k(x,x')$. This is zero since $x\in\uP_\z$ and 
$x'\in\un{U}_P$. Thus $t_{0!}(\cl^{\k'}|_{\Xi_0})=t_{0!}(\bbq)=\p'_{\z\i!}(\bbq)[?]$ and (a) follows.

\subhead 3.2\endsubhead
We preserve the setup of 3.1 and assume that $\vt$ is inner. Let $\z\in\kk^*$. A $\vt$-stable parabolic subgroup 
$P$ of $G$ is said to be {\it $\z$-tight} if $\uP_\z=\un{U}_{P,\z}$. In this case, $P$ is also $\z\i$-tight; 
indeed, we have 
$$\align&\un{U}_{P,\z\i}=\uP^\pe\cap\fg_{\z\i}=\{x\in\fg_{\z\i};\k(x,\uP_\z)=0\}\\&=
\{x\in\fg_{\z\i};\k(x,\un{U}_{P,\z})=0\}=\un{U}_P^\pe\cap\fg_{\z\i}=\uP_{\z\i}.\endalign$$
Now let 

(a) $X$ be a $K$-orbit on the variety of $\vt$-stable parabolic subgroups of $G$ such that some (or equivalently,
any) $P\in X$ is $\z$-tight.
\nl
In this case we have $\tX_\z=\tX'_\z$, $\tX_{\z\i}=\tX'_{\z\i}$ hence $\p_{\z!}(\bbq)=\p'_{\z!}(\bbq)$, 
$\p_{\z\i!}(\bbq)=\p'_{\z\i!}(\bbq)$. Combining this with 3.1(a) we obtain 
$$\fF(\p_{\z!}(\bbq))=\fF(\p'_{\z!}(\bbq))=\p'_{\z\i!}(\bbq)[?]=\p_{\z\i!}(\bbq)[?].\tag b$$
Let $\cb$ be the variety of Borel subgroups of $G$. We show:

(c) {\it if $B\in\cb$, $\vt(B)=B$ and $\z\in\kk^*-\{1\}$ then $B$ is $\z$-tight.}
\nl
Let $\fn=\un{U}_B$. Now $\vt$ acts naturally on $\uB,\fn,\uB/\fn$ and we have an obvious direct sum decomposition
$\uB/\fn=\op_{\z'\in\kk^*}(\uB/\fn)_{\z'}$; moreover, $(\uB/\fn)_{\z'}$ is the image of $\uB_{\z'}$ under 
$\uB@>>>\uB/\fn$. We can find a semisimple element $g_0\in G$ such that $\vt=\Ad(g_0)$. Since $\vt(B)=B$ we have 
$g_0\in B$. Hence the centralizer of $g_0$ in $B$ contains a maximal torus of $B$. Hence $\uB_1$ contains a Cartan
subalgebra $\fh$ of $\uB$. The image of $\fh$ under $\uB@>>>\uB/\fn$ is on the one hand equal to $\uB/\fn$ and on
the other hand is contained in $(\uB/\fn)_1$. Thus $\uB/\fn=(\uB/\fn)_1$. It follows that $(\uB/\fn)_{\z'}=0$ for
any $\z'\ne1$. In particular $(\uB/\fn)_\z=0$. Hence the image of $\uB_\z$ under $\uB@>>>\uB/\fn$ is $0$. In other
words, $\uB_\z\sub\fn$ and (c) follows.

From (c) we see that any $K$-orbit $X$ on $\cb^\vt$ (the variety of $\vt$-stable Borel subgroups of $G$) is as in
(a) hence (b) is applicable to it.

\subhead 3.3\endsubhead
We preserve the setup of 3.1 and assume that $\vt$ is inner. Let $\z\in\kk^*-\{1\}$. We define a collection 
$\fQ'_\z$ of simple perverse sheaves on $\fg_\z$ as follows. A simple perverse sheaf $A$ on $\fg_\z$ is said to be
in $\fQ'_\z$ if there exists a $K$-orbit $X$ on $\cb^\vt$ such that some shift of $A$ is a direct summand of 
$\p_{\z!}\bbq$ where $\p_\z:\tX_\z@>>>\fg_\z$ is defined in terms of $X$ as in 3.1. Note that any object $A$ of 
$\fQ'_\z$ is $K$-equivariant and $\supp(A)$ is contained in $\fg_\z^{nil}$, the variety of nilpotent elements of 
$\fg_\z$ (we use 3.2(b)). 

Note that if $X$ is as in 3.2(a) and $\p_\z:\tX_\z@>>>\fg_\z$ is defined in terms of $X$ as in 3.1 then by the 
decomposition theorem, $\p_{\z!}\bbq\cong\op_hA_h[d_h]$ where $A_h$ are simple perverse sheaves on $\fg_\z$ and 
$d_h\in\ZZ$. We show that 

(a) $A_h\in\fQ'_\z$ for any $h$.
\nl
Let $X'$ be the variety of all $B\in\cb^\vt$ such that $B$ is contained in some (necessarily unique) $P\in X$. 
Define $j:X'@>>>X$ by $B\m P$. Note that $X'$ is a union of (finitely many) $K$-orbits $X'_1,\do,X'_s$ on 
$\cb^\vt$. Let $Y_i=\{(x,B);B\in X'_i;x\in\uB_\z\}$; define $\r_i:Y_i@>>>\fg_\z$ by $(x,B)\m x$. Let 
$Y=\{(x,B);B\in X';x\in\uB_\z\}$; define $\r:Y@>>>\fg_\z$ by $(x,B)\m x$. We have $Y=\sqc_iY_i$ hence 
$\r_!\bbq=\op_i\r_{i!}\bbq$. Define $\s:Y@>>>\tX_\z$ by $\s(x,B)=(x,j(B))$. For $(x,P)\in\tX_\z$ we can identify 
$\s\i(x,P)$ with $\{B\in\cb^\vt;B\sub P\}$ (if $B\in\cb^\vt$, $B\sub P$ we have automatically $x\in\uB_\z$; indeed
we have $x\in\uP_\z=\un{U}_{P,\z}\sub \un{U}_{B,\z}\sub\uB_\z$). We see that $\s$ is a locally trivial fibration 
whose fibres are finite unions of flag manifolds hence $\s_!\bbq\cong\op_j\bbq[2c_j]$ where $c_j$ are integers. We
have $\r=\p_\z\s$ hence $\r_!\bbq=\op_j\p_{\z!}\bbq[2c_j]$. Hence $A_h[?]$ is a direct summand of $\r_!\bbq$. Hence
$A_h[?]$ is a direct summand of $\r_{i!}\bbq$ for some $i$. Thus (a) holds.

Let $A\in\fQ'_\z$. It is known \cite{\VI,\S2, Prop.2} that $\fg_\z^{nil}$ is
a union of finitely many $K$-orbits. It follows that $A=IC(\bco,\ce)_{\fg_\z}[?]$ where $\co$ is a $K$-orbit in 
$\fg_\z^{nil}$, $\bco$ is the closure of $\co$ and $\ce$ is an irreducible $K$-equivariant local system on
$\co$. Thus $A$ is an orbital complex on $\fg_\z$ with unipotent support.

From 3.2(b) we see that $\fF$ defines a bijection from $\fQ'_\z$ (up to isomorphism) to $\fQ'_{\z\i}$ (up to 
isomorphism). We see that 

(b) {\it any object of $\fQ'_\z$ is biorbital. Moreover $\fQ'_\z\sub\fQ_\z$ (see 0.6).}
\nl
We show:

(c) {\it $\fQ'_\z$ is nonempty. Hence $\fQ_\z$ is nonempty. More precisely, there exists a $K$-orbit $\co$ in 
$\fg_\z^{nil}$ such that $IC(\bco,\bbq)_{\fg_\z}[?]\in\fQ'_\z$ where $\bbq$ is viewed as a local system on $\co$.}
\nl
It is well known that $\cb^\vt\ne\em$. Let $X$ be a $K$-orbit on $\cb^\vt$. Define $\p_\z:\tX_\z@>>>\fg_\z$ as in
3.1 in terms of $X$. Note that $\tX_\z$ is vector bundle over $X$ (whose fibre over $B\in X$ is $\uB_\z$); in 
particular $\tX_\z$ is smooth irreducible. Hence $\p_\z(\tX_\z)$ is an irreducible subvariety of $\fg_\z^{nil}$ 
(we use 3.2(b)). There is a unique $K$-orbit $\co$ in $\fg_\z^{nil}$ such that $\co$ is open in $\p_\z(\tX_\z)$. 
Let $n=\dim\p_\z\i(x)$ for any $x\in\co$ (note that $\p_\z\i(x)\ne\em$ for $x\in\co$). For
$x\in\co$ let $S_x$ be the set of irreducible components of dimension $n$ of $\p_\z\i(x)$. We have a finite 
covering $\t:S@>>>\co$ whose fibre at $x\in\co$ is $S_x$. Note that the $2\d$-th cohomology sheaf $\cf$ of 
$(\p_{\z!}\bbq)|_\co$ may be identified with $\t_!\bbq$ (we ignore Tate twists). Hence it contains $\bbq$ as a 
direct summand. Since $(\p_{\z!}\bbq)|_\co$ is a direct sum of shifts of irreducible local systems on $\co$ it 
follows that some shift of $\bbq$ is a direct summand of $(\p_{\z!}\bbq)|_\co$. By the decomposition theorem,
$\p_{\z!}\bbq$ is a direct sum of shifts of irreducible perverse sheaves with support contained in the closure 
$\bco$ of $\co$. Hence some shift of $IC(\bco,\bbq)_{\fg_\z}$ is a direct summand of $\p_{\z!}\bbq$. This proves 
(c).

From (b), (c) we deduce:

(d) {\it there exists a $K$-orbit $\co$ in $\fg_\z^{nil}$ such that $IC(\bco,\bbq)_{\fg_\z}[?]$ belongs to
$\fQ_\z$. Here $\bbq$ is viewed as a local system on $\co$.}
\nl
Let $\VV'_\z$ be the $\QQ(v)$-vector space with basis given by the isomorphism classes of objects in $\fQ'_\z$.
($v$ is an indeterminate). For any $K$-orbit $X$ on $\cb^\vt$ we set $[X]=\sum_{A,j}n_{A,j}v^j\in\VV'_\z$ (sum 
over $A\in\fQ'_\z$ up to isomorphism and $j\in\ZZ$) where $\p_{\z!}\bbq\cong\op_{A,j}A[j]$ 
($\p_\z:\tX_\z@>>>\fg_\z$ is defined in terms of $X$ as in 3.1). We conjecture that:

(e) {\it the elements $[X]$ (for various $K$-orbits $X$ on $\cb^\vt$) generate the vector space $\VV'_\z$.}
\nl
This is known to be true in the case arising from a cyclic quiver (see \cite{\LI}) and also in the case where 
$\z$ has large order in $\kk^*$ (see \cite{\LV}).

\subhead 3.4\endsubhead
We preserve the setup of 0.3 and assume that $\vt$ is inner. Let $\z\in\kk^*-\{1\}$. The following result is
analogous to 3.3(d).

(a) {\it There exists a nonzero function $f:\fg_\z(\FF_q)@>>>\bbq$ (which is constant on each orbit of $K(\FF_q)$)
such that $f$ vanishes on any non-nilpotent element of $\fg_\z(\FF_q)$ and $\hf:\fg_{\z\i}(\FF_q)@>>>\bbq$ 
vanishes on any non-nilpotent element of $\fg_{\z\i}(\FF_q)$.}
\nl
From the assumptions in 0.3, $\cb^\vt$ has at least one $\FF_q$-rational point. Hence it has some irreducible
component $X$ which is defined over $\FF_q$. We have $X(\FF_q)\ne\em$. The morphism $\p_\z:\tX_\z@>>>\fg_\z$ in 
3.1 restricts to a map $t:\tX_\z(\FF_q)@>>>\fg_\z(\FF_q)$. Define $f_\z:\fg_\z(\FF_q)@>>>\bbq$ by 
$f_\z(x)=\sha(t\i(x))$. We have $f_\z(0)=\sha(X(\FF_q))\ne0$. Thus $f_\z\ne0$. From 3.2(b) we deduce 
$\hf_\z=cf_{\z\i}$ where $c\in\bbq^*$ and that $f_\z$ vanishes on any non-nilpotent element of $\fg_\z(\FF_q)$;
similarly, $f_{\z\i}$ vanishes on any non-nilpotent element of $\fg_{\z\i}(\FF_q)$. Hence $\hf_\z$ vanishes on any
non-nilpotent element of $\fg_{\z\i}(\FF_q)$. Clearly, $f_\z$ is constant on each orbit of $K(\FF_q)$. This proves
(a).

\subhead 3.5\endsubhead
Assume that $2\ne0$ in $\kk$. Let $V$ be a $\kk$-vector space of finite dimension $N=2n\ge4$ with a given 
nondegenerate symmetric bilinear form $(,):V\T V@>>>\kk$. Let $K$ be the corresponding special orthogonal group 
acting on $V$ in
an obvious way. For any $\l\in\kk$ let $Q_\l=\{x\in V;(x,x)/2=\l\}$. Note that if $\l\ne0$, then $Q_\l$ is a 
single $K$-orbit in $V$; moreover $Q_0-\{0\}$ is a single $K$-orbit in $V$. Also the isotropy group in $K$ of any
point in $V$ is connected. For $\l\in\kk^*$ let $\ck_\l=IC(Q_\l,\bbq)_V\in\cd(V)$. We set 
$\ck_0=IC(Q_0,\bbq)_V\in\cd(V)$ where $\bbq$ is viewed as a local system on $Q_0-\{0\}$. We set 
$\ck'_0=IC(\{0\},\bbq)_V\in\cd(V)$. Let $V_*=V-Q_0$, an open subset of $V$. For $\l\in\kk^*$ we define 
$\a_\l:V_*@>>>\kk^*$ by $\a_\l(x)=\l(x,x)/2$. We identify $V$ with its dual via $(,)$. Hence 
$\fF:\cd(V)@>>>\cd(V)$ is well defined. The following result describes the antiorbital complexes on $V$.

(i) $\fF(\ck_\l)=IC(V,\a_\l^*\fK^2)[?]$ for any $\l\in\kk^*$;

(ii) $\fF(\ck_0)=\ck_0[?]$;

(iii) $\fF(\ck'_0)=\bbq[?]$.
\nl
Now (iii) is obvious and (ii) is proved in \cite{\LIV}. We prove (i). It is enough to check this at the level
of functions on the set of rational points of $V$ over a finite field. We set $k=\FF_q$. Let $U$ be a $k$-vector
space of dimension $N=2n\ge4$ with a fixed nondegenerate symmetric bilinear form $(,):U\T U@>>>k$ which is split 
over $k$. Let $\l\in k^*$. Define $f:U@>>>\bbq$ by $f(x)=1$ if $(x,x)/2=\l$, $f(x)=0$ if $(x,x)/2\ne\l$. By 0.1,
$\hf:U@>>>\bbq$ is given by $\hf(x)=q^{-n}\sum_{y\in U;(y,y)=2\l}\ps(x,y)$ for $x\in U$. We compute $\hf(x)$ 
assuming that $(x,x)/2=\l'\ne0$. We can find a $2$-dimensional subspace $P$ of $U$ such that $x\in P$ and $(,)$ is
nondegenerate, split on $P$. Let $P'=\{z\in U;(z,P)=0\}$. Note that $(,)$ is nondegenerate split on $P'$. We have 
$$\align&\hf(x)=q^{-n}\sum_{p\in P,p'\in P';(p,p)+(p',p')=2\l}\ps(x,p)\\&=
q^{-n}\sum_{p\in P}\sha(p'\in P';(p',p')=2\l-(p,p))\ps(x,p)\\&=
q^{-n}(\sum_{p\in P;(p,p)=2\l}((q^n-1)(q^{n-1}+1)+1)\ps(x,p)\\&
+\sum_{p\in P;(p,p)\ne2\l}q^{n-1}(q^n-1)\ps(x,p))\\&=
\sum_{p\in P;(p,p)=2\l}\ps(x,p)+q\i(q^n-1)\sum_{p\in P}\ps(x,p).\endalign$$
The last sum is $0$ since $p\m\ps(x,p)$ is a non-trivial character $P@>>>\bbq^*$. Thus
$\hf(x)=\sum_{p\in P;(p,p)=2\l}\ps(x,p)$. We pick a basis $e_1,e_2$ of $P$ such that $(e_1,e_2)=1$, $(e_i,e_i)=0$ 
for $i=1,2$. We set $x=x_1e_1+x_2e_2$ where $x_i\in k$, $x_1x_2=\l'$. We set $p=p_1e_1+p_2e_2$ where $p_i\in k$,
$p_1p_2=\l$. We have
$$\hf(x)=\sum_{p_1,p_2\in k^*;p_1p_2=\l}\ps(x_1p_1+x_2p_2)=\sum_{p'_1,p'_2\in k^*;p'_1p'_2=\l\l'}\ps(p'_1+p'_2).$$
This identity (and the analogous identities where $q$ is replaced by a power of $q$) implies (i).

We now give an alternative proof of (ii) at the level of functions on $U$. For any $n$-dimensional isotropic
subspace $L$ of $U$ we consider the function $f_L:U@>>>\bbq$ which takes the constant value $1$ on $L$ and is $0$
on $U-L$. It is clear that $\hf_L=f_L$. Let $f=c\i\sum_Lf_L$ where $L$ runs over all $n$-dimensional isotropic
subspaces of $U$ and $c=2(q+1)(q^2+1)\do(q^{n-2}+1)$. Clearly, $f(x)=0$ if $(x,x)\ne0$, $f(x)=1$ if $(x,x)=0$, 
$x\ne0$, $f(x)=1+q^{n-1}$ if $x=0$. We have $\hf=f$. This gives the required identity.

Let $\ck$ be a $K$-equivariant simple perverse sheaf on $V$. From (i),(ii),(iii) above we see that the following 
three conditions on $\ck$ are equivalent:

(I) both $\supp(\ck)$ and $\supp(\fF(\ck))$ are contained in $\{x\in V;(x,x)=0\}$;

(II) $\ck$ is biorbital;

(III) $\ck\cong\ck_0[?]$.

\subhead 3.6\endsubhead
Let $V$ be a $\kk$-vector space of dimension $2n+2\ge6$ with a fixed nondegenerate symplectic form 
$\la,\ra:V\T V@>>>\kk$ and with a fixed grading $V=V_0\op V_1$ such that $\la V_0,V_1\ra=0$ and such that
$\dim V_0=2,\dim V_1=2n$. Let $\fs\fp(V)=\{T\in\End(V);\la T(x),y\ra+\la x,T(y)\ra=0\qua\frl x,y\in V\}$. Let
$$E=\{T\in\fs\fp(V);TV_0\sub V_1, TV_1\sub V_0\}.$$
Note that $E=\fs\fp_{-1}$ where $\fs\fp_{\pm1}$ are the $\pm1$ eigenspaces of an involution of $\fs\fp(V)$ 
(induced by an involution of $Sp(V)$ whose fixed point set is $K=Sp(V_0)\T Sp(V_1)$ which acts naturally on $E$).
Hence the notion of antiorbital complex on $E$ is well defined (a special case of 0.6). (In our case the function
$\k_\z$ is the 
restriction to $E\T E$ of the symmetric bilinear form $T,T'\m\tr(TT')$ on $\End(V)$.) The variety of nilpotent 
elements in $E$ decomposes into a union of three $K$-orbits $\{0\},\co,\co'$ represented by $0,N,N'$ where 
$(V,N)\cong(V_{1,3},T_1)^{\op2}\op(V_{1,1},T_1)^{\op(2n-4)}$,
$(V,N')\cong (V_{0,1},T_1)\op(V_{1,0},T_1)\op(V_{1,1},T_1)^{\op(2n-2)}$ (as objects of $\cc^1$ that is without a 
symplectic form). Note that the isotropy groups of $0,N,N'$ in $K$ are connected. Let $\ck_0=IC(\{0\},\bbq)_E[?]$,
$\ck=IC(\bar\co,\bbq)_E[?]$, $\ck=IC(\bar\co',\bbq)_E[?]$ be the corresponding orbital complexes on $E$. We have
the following result.

(a) {\it $\fF(\ck_0)=IC(E,\bbq)[?]$, $\fF(\ck)\cong\ck$, $\fF(\ck')\cong\ck'$. In particular, $\ck$ and $\ck'$ are
biorbital.}
\nl
The first equality in (a) is obvious. Now
let $X$ be the set of all $(T,W)$ where $T\in E$ and $W$ is an $n$-dimensional isotropic subspace of $V_1$ such 
that $T(V_0)\sub W\sub\ker(T)$. Let $X'$ be the set of all $(T,W')$ where $T\in E$ and $W'$ is a line of $V_0$ 
such that $T(V_1)\sub W'\sub\ker(T)$. Note that $X,X'$ are smooth varieties and the obvious projections
$\r:X@>>>E$, $\r':X'@>>>E$ are proper maps.
If $T\in\r(X)$ then clearly $T^3=0$ so that $T$ is nilpotent. Similarly if $T\in\r'(X')$ then $T^3=0$ so that $T$
is nilpotent; actually in this case we have $T^2=0$. (It is enough to show that $T^2V_0=0$ or that
$T^2x=0,T^2x'=0$ where $x,x'$ is a basis of $V_0$ such that $\la x,x'\ra=1$. Clearly, $Tx=0$, $T^2x'=cx$ where
$c\in\kk$. We have $c=\la cx,x'\ra=\la T^2x',x'\ra=-\la Tx',Tx'\ra=0$. Thus $c=0$ and $T^2x'=0$.) Note that 
$\r\i(N)$ can be identified with the variety of all $n$-dimensional isotropic subspaces $W$ of $V_1$ such that $W$
contains $N(V_0)$ (a $2$-dimensional isotropic subspace of $V_1$) and is contained in $\ker(N)\cap V_1$ (a 
codimension $2$ subspace of $V_1$ equal to the perpendicular of $N(V_0)$). Thus $\r$ restricts to a map
$\r\i(\co)@>>>\co$ which is a (locally trivial) fibre bundle whose fibre is isomorphic to the space of Lagrangian 
subspaces of a $(2n-4)$ dimensional symplectic vector space. Since $\co$ is open in the nilpotent variety of $E$ 
we see (using the decomposition theorem) that $\r_!\bbq$ is isomorphic to a direct sum of complexes of the form 
$\ck[?]$ (at least one) and of some complexes of the form $\ck'[?]$ or $\ck_0[?]$.

Next we note that $\r'{}\i(N')$ is a single point of $X'$ namely $(N',W')$ where $W'=N'(V_1)$; more precisely, the
restriction of $\r'$ from $\r'{}\i(\co')$ to $\co'$ is an isomorphism. Since $\r'(X')\sub\co'\cup\{0\}$ and $\co'$
is open in $\co'\cup\{0\}$ we see (using the decomposition theorem) that $\r'_!\bbq$ is isomorphic to $\ck'[?]$ 
direct sum with some complexes of the form $\ck_0[?]$.

Let $W$ be an $n$-dimensional isotropic subspace of $V_1$ and let $W'$ be a line in $V_0$. Let 

$\fp=\{T\in\fs\fp(V);TW\sub W\}$, $\fp'=\{T\in\fs\fp(V);TW'\sub W'\}$,

$\fn=\{T\in\fs\fp(V);TV\sub W^\pe,TW^\pe\sub W,TW=0\}$, 
$\fn'=\{T\in\fs\fp(V);TV\sub W'{}^\pe, TW'{}^\pe\sub W',TW'=0\}$. 
\nl
Here $W^\pe,W'{}^\pe$ denote the perpendicular 
to $W,W'$ with respect to $\la,\ra$. Note that $\fp,\fp'$ are parabolic subalgebras of $\fs\fp(V)$ with 
nil-radicals $\fn,\fn'$. Moreover we have $\fp=\fp_1\op\fp_{-1},\fp'=\fp'_1\op\fp'_{-1}$, $\fn=\fn_1\op\fn_{-1}$,
$\fn'=\fn'_1\op\fn'_{-1}$ where $()_i=()\cap\fs\fp_i$. We show that $\fp_{-1}=\fn_{-1}$, $\fp'_{-1}=\fn'_{-1}$
(that is, the parabolic subgroups of $Sp(V)$ corresponding to $\fp,\fp'$ are $(-1)$tight). 
Let $T\in\fp_{-1}$. We have $TW\sub V_0\cap W=0$ so that $TW=0$. It follows that $TV\sub W^\pe$. We have 
$W^\pe=V_0\op W$ hence $TW^\pe=TV_0+TW=TV_0\sub W^\pe\cap V_1=W$. We see that $\fp_{-1}=\fn_{-1}$. A similar proof
shows that $\fp'_{-1}=\fn'_{-1}$. Now $\r_!\bbq=\tInd_{\fp_{-1}}^{\fs\fp_{-1}}\bbq$,
$\r'_!\bbq=\tInd_{\fp'_{-1}}^{\fs\fp_{-1}}\bbq$ (see A.1). Here $\bbq$ is viewed as a complex on
$\fp_{-1}/\fn_{-1}=0$ or $\fp'_{-1}/\fn'_{-1}=0$. From Theorem A.2 (or from 3.2(b)) we see that 
$\fF(\r_!\bbq)=\r_!\bbq[?]$, $\fF(\r'_!\bbq)=\r'_!\bbq[?]$. 

It follows that $\fF(\r'_!\bbq)$ is isomorphic to $\ck'[?]$ direct sum with some complexes of the form $\ck_0[?]$;
it is also isomorphic to $\fF(\ck')[?]$ direct sum with some complexes of the form $\fF(\ck_0)[?]$. But 
$\supp\fF(\ck_0)=E$ showing that $\fF(\ck_0)[?]$ cannot be a direct summand of $\fF(\r'_!\bbq)$. It follows that 
$\r'_!\bbq\cong\ck'[?]\cong\fF(\ck')[?]$. Hence $\fF(\ck')\cong\ck'$.

We also see that $\fF(\r_!\bbq)$ is isomorphic to a direct sum of complexes $\ck[?]$ (at least one) and complexes
of the form $\ck'[?]$ or $\ck_0[?]$; it is also isomorphic to a direct sum of complexes $\fF(\ck)[?]$ (at least 
one) and complexes of the form $\fF(\ck')[?]=\ck'[?]$ or $\fF(\ck_0)[?]$. Again $\fF(\ck_0)[?]$ cannot be a direct
summand of $\fF(\r_!\bbq)$. It follows that $\fF(\ck)$ is isomorphic to $\ck$ or to $\ck'$. If $\fF(\ck)\cong\ck'$
then $\ck\cong\fF(\ck')$ hence $\ck\cong\ck'$ which is not the case. Hence we have $\fF(\ck)\cong\ck$. This proves
(a).

\subhead 3.7\endsubhead
Let $V$ be a $\kk$-vector space of dimension $2n$ with a fixed nondegenerate symplectic form 
$\la,\ra:V\T V@>>>\kk$. Let 
$$E=\{T\in\End(V);\la T(x),y\ra=\la x,T(y)\ra\qua\frl x,y\in V\}.$$
Note that $E$ can be viewed as the $(-1)$ eigenspace of an involution of $\End(V)$ (induced by an involution of
$GL(V)$ whose fixed point set is the symplectic group $Sp(V)$ which acts naturally on $E$). Hence the notion of 
antiorbital complex on $E$ is well defined (a special case of 0.6). (In our case the function $\k_\z$ is the
restriction to $E\T E$ of the symmetric bilinear form $T,T'\m\tr(TT')$ on $\End(V)$.) Let $E_0$ be the set of all
$T\in E$ such that $T:V@>>>V$ is semisimple and any eigenspace of $T$ is $2$-dimensional. Note that $E_0$ is open
dense in $E$. Using methods similar to those in \S2 we see that any antiorbital complex on $E$ is of the form
$IC(E_0,\cl)[?]$ for a suitable local system $\cl$ on $E_0$ (compare with \cite{\LIV, \S13}.) It follows that, if
$n>0$, there are no biorbital complexes on $E$.

\subhead 3.8\endsubhead
We preserve the setup of 3.1. Let $\z\in\kk^*$. For any $K$-orbit $\co$ in $\fg_\z^{nil}$ let
$\co^!=\{(x,y)\in\fg_\z\T\fg_{\z\i};x\in\co,[x,y]=0\}$ where $[,]$ is the bracket in $\fg$. As in 
\cite{\LV, 22.2}, we identify $\co^!$ with the conormal bundle of $\co$ in $\fg_\z$. Hence it is smooth,
irreducible of dimension $\dim\fg_\z$. Hence $A:=\{(x,y)\in\fg_\z^{nil}\T\fg_{\z\i};[x,y]=0\}$ is a (closed) 
subvariety of $\fg_\z\T\fg_{\z\i}$ of pure dimension $\dim\fg_\z$; its irreducible components are $\ov{\co^!}$ 
(closure of $\co^!$) for various $\co$ as above. Similarly for any $K$-orbit $\cv$ in $\fg_{\z\i}^{nil}$ let
$\cv^!=\{(x,y)\in\fg_\z\T\fg_{\z\i};y\in\cv,[x,y]=0\}$. Then $\cv^!$ is smooth, irreducible of dimension  
$\dim\fg_{\z\i}=\dim\fg_\z$. Hence $A':=\{(x,y)\in\fg_\z\T\fg_{\z\i}^{nil};[x,y]=0\}$ is a (closed) subvariety of 
$\fg_\z\T\fg_{\z\i}$ of pure dimension $\dim\fg_\z$; its irreducible components are $\ov{\cv^!}$ (closure of 
$\cv^!$) for various $\cv$ as above. Let 
$$\L=A\cap A'=\{(x,y)\in\fg_\z^{nil}\T\fg_{\z\i}^{nil};[x,y]=0\},$$
a (closed) subvariety of $\fg_\z\T\fg_{\z\i}$ of dimension $\le\dim\fg_\z$. The irreducible components of $\L$ of
dimension $\dim\fg_\z$ are of the form $\ov{\co^!}$ (where $\co$ runs through a subset $H_\z$ of the set of 
$K$-orbits on $\fg_\z^{nil}$). They are also of the form $\ov{\cv^!}$ (where $\cv$ runs through a subset 
$H_{\z\i}$ of the set of $K$-orbits on $\fg_{\z\i}^{nil}$). Hence there is a unique bijection 
$\io:H_\z@>>>H_{\z\i}$ such that $\ov{\co^!}=\ov{\io(\co)^!}$ for any $\io\in H_\z$.

Now let $\ck\in\fQ_\z$ and assume that $\supp\ck=\bco$ (closure of $\co$ as above), $\supp\fF(\ck)=\bar\cv$ 
(closure of $\cv$ as above). Let $S$ (resp. $S'$) be the singular support of $\ck$ (resp. $\fF(\ck)$); they are 
subvarieties of $\fg_\z\T\fg_{\z\i}$. Note that $S=S'$ and $\co^!\sub S\sub\cup_{\co'\sub\bco}\co'{}^!\sub A$, 
$\cv^!\sub S'\sub\cup_{\cv'\sub\bar\cv}\cv'{}^!\sub A'$ where $\co'$ runs over the $K$-orbits in $\bco$ and $\cv'$
runs over the $K$-orbits in $\bar\cv$. It follows that $\co^!\sub A'$. Since $\co^!\sub A$ we see that 
$\co^!\sub A\cap A'=\L$. Similarly, $\cv^!\sub\L$. We see that:

(a) {\it if $\ck\in\fQ_\z$ then $\supp(\ck)=\bco$ where $\co\in H_\z$. In particular if $x\in\co$ and
$y\in\fg_{\z\i}$ satisfies $[x,y]=0$ then $y$ is nilpotent.}
\nl
In the case where $\fg_\z$ arises as in 0.4 from a cyclic quiver, the variety $\L$ is the same as that defined in
\cite{\LI, \S12}. In this case $\L$ is of pure dimension. Also in this case the bijection $\io$ can be viewed as a
bijection $\cp^{ap}_\nu@>\si>>\cp^{ap}_\nu$ (notation of 1.2). It would be interesting to describe this bijection
explicitly. (It is analogous to the involution of $\car_t$ which takes a partition to the conjugate partition.)

If $\z$ has large order in $\kk^*$ then $\L$ is again of pure dimension (see \cite{\LV, 22.2}).

\subhead 3.9\endsubhead
Let $G$ (resp. $\fg$) be the set of all $4\T4$ matrices $a=(a_{ij})$ with $a_{ij}\in\kk$ for $i,j\in[1,4]$,
$a_{ij}=0$ for all $i>j$ and $a_{ii}=1$ (resp. $a_{ii}=0$) for all $i$. Note that $G$ is naturally a unipotent
algebraic group with Lie algebra $\fg$. Let $\fh$ be the set all $4\T4$ matrices $b=(b_{ij})$ with 
$b_{ij}\in\kk$ for $i,j\in[1,4]$, $b_{ij}=0$ for all $i\le j$. We identify $\fg$ with the dual space of $\fh$ via
the nondegenerate bilinear pairing $(a,b)=\sum_{i<j}a_{ij}b_{ji}\in\kk$. Now $G$ acts on $\fg$ by conjugation and
this induces a $G$-action on $\fh$ (by duality). Hence the orbital complexes on $\fh$ are well defined and the 
antiorbital complexes on $\fg$ are well defined.

Consider the partition $\fh=U_1\cup U_2\cup U_3\cup U_4\cup U_5$ where 

$U_1=\{b;b_{31}=b_{41}=b_{42}=0\}$; $U_2=\{b;b_{31}=b_{41}=0, b_{42}\ne0\}$;

$U_3=\{b;b_{41}=b_{42}=0, b_{31}\ne0\}$; $U_4=\{b;b_{41}\ne0\}$; 

$U_5=\{b;b_{41}=0,b_{31}\ne0,b_{42}\ne0$.
\nl
For $i=1,2,3,4,5$, $U_i$ is a union of $G_n$-orbits of fixed dimension: $0,2,2,4,2$.

Let 

$V_1=\fg$, $V_2=\{a\in\fg;a_{23}=a_{34}=0\}$, $V_3=\{a\in\fg;a_{12}=0,a_{23}=0\}$, 

$V_4=\{a\in\fg;a_{12}=a_{34}=0\}$, $V'_4=\{a\in V_4; a_{23}\ne0\}$, 

$V_{5;x,y}=\{a\in\fg;a_{23}=0,xa_{12}-ya_{34}=0\}$ ($x\in\kk^*,y\in\kk^*$).
\nl
The class of antiorbital complexes on $\fg$ consists of: 

$(1)$ $\cl^f[6]$ where $f:V_1@>>>\kk$ is $a\m xa_{12}+ya_{23}+za_{34}$ ($x,y,z\in\kk$);

$(2)$ $\cl^f_\fg[4]$ where $f:V_2@>>>\kk$ is $a\m xa_{12}+ya_{24}$, ($x\in\kk,y\in\kk^*$);

$(3)$ $\cl^f_\fg[4]$ where $f:V_3@>>>\kk$ is $a\m xa_{13}+ya_{34}$, ($x\in\kk^*,y\in\kk$);

$(4)$ $IC(V_4,\cl^f)_\fg[4]$ where $f:V'_4@>>>\kk$ is $a\m xa_{23}+y(a_{14}-\fra{a_{24}a_{13}}{a_{23}})$, 
($x\in\kk,y\in\kk^*$);

$(5)$ $\cl^f_\fg[4]$ where $f:V_{5,c,d}@>>>\kk$ is $a\m xa_{13}+ya_{24}+zx\i a_{34}$, ($z\in\kk$).
\nl
Note that the complexes $(1),(2),(3),(4),(5)$ are obtained by applying $\fF$ to the orbital complexes on $\fh$ 
with support contained in $U_1,U_2,U_3,U_4,U_5$ respectively.

\head Appendix. Induction, restriction\endhead
\subhead A.1\endsubhead
Let $G,\fg,\vt,K,\k$ be as in 0.3. We assume that the characteristic of $\kk$ is sufficiently large. For 
$\z\in\kk^*$ let $\fg_\z,\k_\z$ be as in 0.3. We consider the datum 

$(P,U,l,\fp,\fn,\fl)$
\nl
where $P$ is a parabolic subgroup of $G$ such that $\vt(P)=P$, $U$ is the unipotent radical of $P$, $L=P/U$ and 
$\fp,\fn,\fl$ denote the Lie algebras of $P,U,L$. Note that $\fp=\op_{\z\in\kk^*}\fp_\z$,
$\fn=\op_{\z\in\kk^*}\fn_\z$ where $\fp_\z=\fp\cap\fg_\z$, $\fn_\z=\fn\cap\fg_\z$. Moreover $\vt$ induces a 
semisimple automorphism of $L$ and of $\fl$ denoted again by $\vt$ and we have a decomposition 
$\fl=\op_{\z\in\kk^*}\fl_\z$ where $\fl_\z$ is the $\z$-eigenspace of $\vt:\fl@>>>\fl$ (it is also equal to
$\p(\fp_z)$ where $\p:\fp@>>>\fl$ is the obvious map). Let $L_K$ be the identity component of the fixed point set 
of $\vt:L@>>>L$. It acts naturally on $\fl_\z$ for any $\z\in\kk^*$.

Let $\z\in\kk^*$. Let $\p_\z:\fp_\z@>>>\fl_\z$ be the restriction of $\p$. Let $U_K=U\cap K$ (a connected 
unipotent group), $P_K=(P\cap K)^0$. Let $E'_\z=K\T_{U_K}\fp_\z,E''_\z=K\T_{P_K}\fp_\z$. we have a diagram
$$\fl_\z@<p_1^\z<<E'_\z@>p_2^\z>>E''_\z@>p_3^\z>>\fg_\z$$
where $p_i=p_i\z$ are given by $p_1(g,x)=\p_\z(x)$, $p_2(g,x)=(g,x)$, $p_3(g,x)=\Ad(g)(x)$. Note that $p_1$ is 
smooth with connected fibres, $p_2$ is a principal bundle with group $P_K/U_K=L_K$ and $p_3$ is a proper map. Let
$A$ be a semisimple $L_K$-equivariant complex on $\fl_\z$. Then $p_1^*A$ is a semisimple $L_K$-equivariant complex
on $E'_\z$. Hence $p_1^*A=p_2^*A'$ for a well defined semisimple complex $A'$ on $E''_\z$. We set
$$\Ind_{\fp_\z}^{\fg_\z}(A)=p_{3!}A'\in\cd(\fg_\z),$$
$$\tInd_{\fp_\z}^{\fg_\z}(A)=\Ind_{\fp_\z}^{\fg_\z}(A)[\dim\fn_0+\dim\fn_\z].$$
The following result is of the same type as \cite{\LIII, 7(a)}, \cite{\LI, Theorem 5.4}, \cite{\LV, Cor.10.5},
\cite{\HE, Theorem 4.1}. The proof we give is almost a word by word repetition of that of Theorem 5.4 in 
\cite{\LI}. 

\proclaim{Theorem A.2} We preserve the setup of A.1. We have
$$\fF(\tInd_{\fp_\z}^{\fg_\z}(A))\cong\tInd_{\fp_{\z\i}}^{\fg_{\z\i}}(\fF(A)).$$
\endproclaim
We consider the  commutative  diagram
$$\CD
X_a     @<u_{ba}<<      X_b  @>u_{bc}>>        X_c     @>u_{cd}>>       X_d \\
@.                 @Au_{eb}AA               @Au_{fc}AA            @Au_{gd}AA\\
{}   @.                 X_e  @>u_{ef}>>          X_f @>u_{fg}>>         X_g \\
@Au_{ja}AA            @Au_{he}AA                 @Au_{if}AA              @. \\
{}         @.           X_h  @>u_{hi}>>           X_i         @.         {} \\
@.                  @Vu_{hk}VV               @Vu_{il}VV      @Vu_{gp}VV     \\
X_j    @<u_{kj}<<       X_k @>u_{kl}>>              X_l       @.    {}   \\
@Vu_{jm}VV               @Vu_{kn}VV                @Vu_{lo}VV   @.   \\
X_m      @<u_{mn}<<      X_n    @>u_{no}>>         X_o        @>u_{op}>> X_p
\endCD$$
in which the notation is as follows.

$X_a=\fl_\z$, $X_c=E''_\z$, $X_b=E'_\z$, $X_d=\fg_\z$.

$X_m=\fl_{\z\i}$, $X_o=E''_{\z\i}$, $X_n=E'_{\z\i}$, $X_p=\fg_\z$.

$X_g=\fg_\z\T\fg_{\z\i}$, $X_f=K\T_{P_K}(\fp_\z\T\fg_{\z\i})$, $X_e=K\T_{U_K}(\fp_\z\T\fg_{\z\i})$.

$X_i=K\T_{P_K}(\fp_\z\T\fp_{\z\i})$, $X_h=K\T_{U_K}(\fp_\z\T\fp_{\z\i})$.

$X_j=\fl_\z\T\fl_{\z\i}$, $X_k=K\T_{U_K}(\fl_\z\T\fp_{\z\i})$, $X_l=K\T_{P_K}(\fl_\z\T\fp_{\z\i})$.

$u_{ba}=p_1^\z,u_{bc}=p_2^\z,u_{cd}=p_3^\z$, $u_{mn}=p_1^{\z\i},u_{no}=p_2^{\z\i},u_{op}=p_3^{\z\i}$.

$u_{fg}$ is $(g,x,x')\m(\Ad(g)x,Ad(g)x')$, $u_{kj}$ is $(g,x,x')\m(x,\p_{\z\i}(x'))$.

$u_{ef},u_{hi},u_{kl},u_{eb},u_{fc},u_{gd},u_{ja},u_{he},u_{if},u_{hk},u_{il},u_{gp},u_{kn},u_{lo},u_{jm}$ 
are the obvious maps.

Note that $\k:\fg\T\fg@>>>\kk$ restricted to $\fp\T\fp$ induces a function $\fl\T\fl@>>>\kk$ which is analogous to
$\k$. Hence we can define local systems of rank $1$ on $X_j,X_k,X_h,X_l,X_i,X_e,X_f,X_g$ (equal to $\cl^{\k_\z}$ 
in the case of $X_g$) which correspond to each other under inverse image by
$$u_{fg},u_{ef},u_{he},u_{if},u_{hi},u_{kj},u_{hk},u_{kl},u_{il}.$$
We shall denote each of these local systems by $\cl$. 

Let $L_a=A\in\cd(X_a)$. Let $L_b=u_{ba}^*L_a\in\cd(X_b)$; let $L_c\in\cd(X_c)$ be the unique semisimple complex 
such that $u_{bc}^*L_c=L_b$ and let $L_d=u_{cd!}L_c\in\cd(X_d)$. Let $L_g=u_{gd}^*L_d\in\cd(X_g)$ and let
$L_p=u_{gp}!(L_g\ot\cl)\in\cd(X_p)$. By definition we have $\fF(\Ind_{\fp_\z}^{\fg_\z}(A))=L_p[D]$ where 
$D=\dim\fg_\z$. Now let $L_j=u_{ja}^*L_a\in\cd(X_j)$ and $L_m=u_{jm!}(L_j\ot\cl)\in\cd(X_m)$. This is a 
semisimple, $L_K-$equivariant complex (it is $\fF(A)[-\dim\fl_\z]$). Let $L_n=u_{nm}^*L_m\in\cd(X_n)$ and let 
$L_o\in\cd(X_o)$ be the unique semisimple complex such that $u_{no}^*L_o=L_n$. Let $L'_p=u_{op!}L_o\in\cd(X_p)$. 
By definition we have $\Ind_{\fp_{\z\i}}^{\fg_{\z\i}}(\fF(A))=L'_p[\dim\fl_\z]$. Hence it suffices to prove that 
$$L_p\cong L'_p[\dim\fl_\z-\dim\fg_\z+\dim\fn_{\z\i}-\dim\fn_\z].\tag a$$
Let $L_e=u_{eb}^*L_b\in\cd(X_e)$, $L_f=u_{fc}^*L_c\in\cd(X_f)$. Then $L_f$ is a semisimple complex (since $L_c$ is
semisimple and $u_{fc}$ is smooth with connected fibres) and $L_e=u_{ef}^*L_f$. Moreover, $u_{fg!}L_f=L_g$ (since
the diagram $u_{fc},u_{cd},u_{fg},u_{gd}$ is cartesian). Hence we may go from $L_a$ to $L_p$ by the chain 
$L_e=(u_{ba}u_{eb})^*L_a$, $L_e=u_{ef}^*L_f$ ($L_f$ semisimple), $L_p=(u_{gp}u_{fg})_!(L_f\ot\cl)$. Similarly, we
may go from $L_a$ to $L'_p$ by the chain $L_k=(u_{ja}u_{kj})^*L_a\in\cd(X_k)$, $L_k=u_{kl}^*L_l$ ($L_l\in\cd(X_l)$
semisimple), $L'_p=(u_{op}u_{lo})_!(L_l\ot\cl)$.

Let $L_h=u_{hk}^*L_k\in\cd(X_h)$, $L_i=u_{il}^*L_l\in\cd(X_i)$. Note that $u_{il}$ is a vector bundle with fibres
of dimension $\dim\fn_\z$. It follows that $L_i$ is semisimple (recall that $L_l$ is semisimple) and that 
$u_{il!}L_i=L_l[-2\dim\fn_\z]$. Hence we have $(u_{op}u_{lo}u_{il})_!(L_i\ot\cl)=L'_p[-2\dim\fn_\z]$. We have
$\dim\fn_{\z\i}-\dim\fn_\z=\fg_\z-\fl_\z-2\dim\fn_\z$. Hence it is enough to prove that 
$$L_p\cong L'_p[-2\dim\fn_\z].\tag a1$$
We have $u_{hi}^*L_i=L_h$ and it follows that we  may go from $L_a$ to $L'_p$ by the  chain
$L_h=(u_{ja}u_{kj}u_{hk})^*L_a$, $L_h=u_{hi}^*L_i$ ($L_i\in\cd(X_i)$ semisimple), 
$L'_p[-2\dim\fn_\z]=(u_{op}u_{lo}u_{il})_!(L_i\ot\cl)$.

Since $L_c$ is semisimple and $u_{fc}u_{if}$ a vector bundle we see that $u_{if}^*L_f=(u_{fc}u_{if})^*L_c$ is 
semisimple. Now both $u_{if}^*L_f$ and $L_i$ are semisimple and they have the same inverse image $L_h$ under 
$u_{hi}$ (a smooth morphism with connected fibres). It follows that $u_{if}^*L_f\cong L_i$. Since 
$u_{op}u_{lo}u_{il}=u_{gp}u_{fg}u_{if}$ we see that 
$L'_p[-2\dim\fn_\z]=(u_{gp}u_{fg}u_{if})_!(u_{if}^*(L_f\ot\cl))$. We now see that (a1) would be a consequence of 
the following statement:

$L_f\ot\cl$ and $u_{if!}u_{if}^*(L_f\ot\cl)$ have the same image under $(u_{gp}u_{fg})_!$. 

An equivalent statement is: 

(b) if $u'$ denotes the inclusion of $X_f-X_i$ into $X_f$ (as an open subset), then \lb
$(u_{gp}u_{fg}u')_!u'{}^*(L_f\ot\cl)=0$.

(We use the distinguished triangle associated with the partition $X_f=X_i\cup(X_f-X_i)$.) We now consider the 
commutative diagram
$$\CD
X_b  @>u_{bc}>> X_c @<u_{fc}<< X_f @<u'<< X_f-X_i        \\
@Vu_bVV         @Vu_cVV       @Vu_fVV     @Vu''VV   \\
Y_b@>w_{bc}>> Y_c @<w_{fc}<< Y_f@<w'<< Y_f-X_l \\
@Vv'VV             @.               @Vw''VV       @.      \\
X_a           @.             @.                 X_p  @.
\endCD$$
where the notation is as follows.

$Y_b=(K/U_K)\T\fl_\z$, $Y_c=K\T_{P_K}\fl_\z$, $Y_f=K\T_{P_K}(\fl_\z\T\fg_{\z\i})$.

$u_b,u_c,u_f,w_{bc},w_{fc},w',v'$ are the obvious maps, $w''$ is $(g,x,x')\m\Ad(g)x'$.

Let $M_b=v'{}^*L_a\in\cd(Y_b)$; this is a $L_K$-equivariant semisimple complex hence there is a well defined
semisimple complex $M_c\in\cd(Y_c)$ such that $w_{bc}^*M_c=M_b$. Let $M_f=w_{fc}^*M_c\in\cd(Y_f)$. It is clear 
that $u_b^*M_b=L_b$, $u_c^*M_c=L_c$, $u_f^*M_f=L_f$ (note that $u_b,u_c,u_f$ are vector bundles). Hence we have
$u'{}^*L_f=u'{}^*u_f^*M_f=u''{}^*w'{}^*M_f$. The statement (b) can now be rewritten in terms of $M_f$ instead of
$L_f$:

$(u_{gp}u_{fg}u')_!(u''{}^*w'{}^*M_f \ot\cl)=0$

or equivalently (using $u_{gp}u_{fg}u'=w''w'u''$):

$(w''w'u'')_!(u''{}^*w'{}^*M_f \ot\cl)=0$.

This would be a consequence of the following statement:

$u''_!(u''{}^*w'{}^*M_f \ot\cl)=0$.

We have $u''_!(u''{}^*w'{}^*M_f\ot\cl)=w'{}^*M_f \ot(u''_!\cl)$, hence it suffices to prove that
   
$u''_!\cl=0$ in $\cd(Y_f-X_l)$.

Let us fix a point $(g,x,x')\in Y_f-X_l$ and let $\G$ be the fibre of $w'$ over this point. Let $\ti\k:\G\to\kk$ 
be the function $(g,x,x')\m\k(x,x')$. By base change, it is enough to prove that the cohomology with compact 
support of $\G$ with coefficients in $\cl^{\ti\k}|_\G$ is zero. By a known property of Artin-Schreier local 
systems, it is enough to verify the following statement: one can identify $\G$ with $\kk^N$ for some $N$ so that 
$\ti\k$ is given by a non-constant affine linear form on $\kk^N$. We can find a $\vt$-stable Levi subgroup $\tL$ 
of $P$. Let $\ti\fl$ be the Lie algebra of $\tL$. Let $\ti\fl_\z=\ti\fl\cap\fg_\z$. We have 
$\fp_\z=\ti\fl_\z\op\fn_\z$. We can identify $\fl_\z=\ti\fl_\z$ in an obvious way; hence we can view $x$ as an 
element of $\ti\fl_\z$. We have an isomorphism $\fn_\z@>\si>>\G$ given by $z\m(g,x+z,x')$. The function $\ti\k$ on
$\G$ can be identified with the function $z\m\k(x+z,x')=\k(z,x')+\k(x,x')$ on $\fn_\z$. It is enough to show that
the linear function $z\m\k(z,x')$ is not identically zero. Assume that it is identically zero that is 
$\k(\fn_\z,x')=0$. If $\z'\ne\z$ we have $\k(\fg_{\z'},\fg_{\z\i})=0$. Hence $\k(\fn_{\z'},x')=0$. We see that
$\k(\fn,x')=0$ hence $x'\in\fp'$ and $x'\in\fp_{\z\i}$. This contradicts $x'\in\fg_{\z\i}-\fp_{\z\i}$. The theorem
is proved.

\subhead A.3\endsubhead
In the setup of A.1 let us assume that $A$ is an orbital complex (on $\fl_\z$). We show:

(a) {\it $\tInd_{\fp_\z}^{\fg_\z}(A)$ is a direct sum of finitely many orbital complexes (with shifts) on 
$\fg_\z$.}
\nl
Let $A',p_3$ be as in A.1. In our case $A'$ can be viewed as a mixed complex, pure of weight $0$ (up to shift); 
hence by the decomposition theorem, $p_{3!}A'$ is a direct sum of shifts of simple perverse sheaves. It is then
enough to show that $S:=\supp(p_{3!}A)$ is contained in the union of finitely many $K$-orbits on $\fg_\z$. Since 
the intersection of any $G$-orbit on $\fg$ with $\fg_\z$ is a union of finitely many $K$-orbits, it is enough to 
show that $S$ is contained in the union of finitely many $G$-orbits on $\fg$. Let $S_0=\supp(A)\sub\fl_\z$. From 
the definitions, $S$ is contained in the union of $K$-orbits in $\fg_\z$ that meet $\p_\z\i(S_0)$. Hence $S$ is 
contained in the union of $G$-orbits in $\fg$ that meet $\p\i(S_0)$. Hence it is enough to show that $\p\i(S_0)$
is contained in a finite union of $G$-orbits in $\fg$. It is also enough to show that the set $\cs$ of
semisimple parts of the various elements in $\p\i(S_0)$ is contained in a single $G$-orbit in $\fg$. Now the set
$\cs_0$ of semisimple parts of the various elements in $S_0$ is contained in a single $L_K$-orbit in $\fl_\z$ 
hence in a single $L$-orbit in $\co_0$ in $\fl$. Since $\cs\sub\p\i(\cs_0)$ it is enough to note that the set of
semisimple elements in $\p\i(\cs_0)$ is a single $P$-orbit on $\fp$. This completes the proof of (a).

Next we assume, in the setup of A.1, that $A$ is an antiorbital complex (on $\fl_\z$). We show:

(b) {\it $\tInd_{\fp_\z}^{\fg_\z}(A)$ is a direct sum of finitely many antiorbital complexes (with shifts) on 
$\fg_\z$.}
\nl
By assumption we have $A=\fF(A_1)$ where $A_1$ is an orbital complex on $\fl_{\z\i}$. Using A.2 with $A$ replaced
by $A_1$ and $\z$ by $\z\i$ we see that 

$\tInd_{\fp_\z}^{\fg_\z}(A)=\fF(\tInd_{\fp_{\z\i}}^{\fg_{\z\i}}(A_1))$.
\nl
Using (a) with $A$ replaced by $A_1$ and $\z$ by $\z\i$ we see that 

$\tInd_{\fp_{\z\i}}^{\fg_{\z\i}}(A_1)\cong\op_iP_i[d_i]$
\nl
where $P_i$ are orbital complexes on $\fg_{\z\i}$ and $d_i\in\ZZ$. It follows that 

$\tInd_{\fp_\z}^{\fg_\z}(A)\cong\op_i\fF(P_i)[d'_i]$
\nl
where $d'_i\in\ZZ$ and (b) is proved.

Finally we assume, in the setup of A.1, that $A$ is a biorbital complex (on $\fl_\z$). We show:

(c) {\it $\tInd_{\fp_\z}^{\fg_\z}(A)$ is a direct sum of finitely many biorbital complexes (with shifts) on 
$\fg_\z$.}
\nl
Using (a) we see that $\tInd_{\fp_\z}^{\fg_\z}(A)\cong\op_iP_i[d_i]$ where $P_i$ are orbital complexes on 
$\fg_\z$ and $d_i\in\ZZ$. Using (b) we see that $\op_iP_i[d_i]\cong\op_jQ_j[e_j]$ where $Q_j$ are antiorbital 
complexes on $\fg_\z$ and $e_j\in\ZZ$. It follows that the isomorphism classes in $\{P_i\}$ coincide with the 
isomorphism classes in $\{Q_j\}$. Hence each $P_i$ is biorbital. This proves (c).

\subhead A.4\endsubhead
Let $\z\in\kk^*$. For $B\in\cd(\fg_\z)$ we define (in the setup of A.1) 
$\Res_{\fp_\z}^{\fg_\z}(B)=\p_{\z!}(B|_{\fp_\z})\in\cd(\fl_\z)$,
$\tRes_{\fp_\z}^{\fg_\z}(B)=\Res_{\fp_\z}^{\fg_\z}(B)[\dim\fn_\z-\dim\fn_1]$.

\proclaim{Proposition A.5} We have $\fF(\tRes_{\fp_\z}^{\fg_\z}(B))\cong
\tRes_{\fp_{\z\i}}^{\fg_{\z\i}}(\fF(B))$ in $\cd(\fl_{\z\i})$.
\endproclaim
The proof is exactly the same as that in \cite{\LV, 10.4} if we use Lemma A.6 below instead of \cite{\LV, 10.3} 
and we replace $\fg_n,\fp_n,\fl_n,\fn_n,\fg_{-n},\fp_{-n},\fl_{-n},\fn_{-n},\fn_0,\t,\cl_\t$ by
$\fg_\z,\fp_\z,\fl_\z,\fn_\z,\fg_{\z\i},\fp_{\z\i},\fl_{\z\i},\fn_{\z\i},\fn_1,\k_\z,\cl^{\k_\z}$.

\proclaim{Lemma A.6} Consider the diagram
$$\CD
(\fg_\z\T\fp_{\z\i})-(\fp_\z\T\fp_{\z\i})@>>>\fg_\z\T\fp_{\z\i}\\
@V\c'VV                                             @V\c VV   \\
(\fg_\z\T\fl_{\z\i})-(\fp_\z\T\fl_{\z\i})@>>>\fg_\z\T\fl_{\z\i}\endCD$$
where the horizontal maps are the obvious inclusions, $\c$ is the obvious projection and $\c'$ makes the diagram 
commutative. We have $\c'_!(\cl^{\k_\z})=0$. (Here $\cl^{\k_\z}$ is the local system on
$(\fg_\z\T\fp_{\z\i})-(\fp_\z\T\fp_{\z\i})$ defined by restricting the analogous local system on
$\fg_\z\T\fp_{\z\i}$, see 0.3.)
\endproclaim
We choose a $\vt$-stable Levi subgroup $\tL$ of $P$. Let $\ti\fl$ be the Lie algebra of $\tL$. We have
$\ti\fl=\op_{\z'\in\kk^*}\ti\fl_{\z'}$ where $\ti\fl_{\z'}=\ti\fl\cap\fg_{\z'}$. We identify $\ti\fl$ with $\fl$ 
via the obvious map $\fp@>>>\fl$. Let $(x,y)\in\fg_\z\T\fl_{\z\i}$ be such that $x\n\fp_\z$. We have

(a) $\c'{}\i(x,y)=\{(x,y+z);z\in\fn_{\z\i}\}$.
\nl
It is enough to show that the cohomology with compact support of the variety (a) with coefficients in 
$\cl^{\k_\z}$ is zero. We identify this variety with $\fn_{\z\i}$ via the coordinate $z$. The restriction of
$\k_\z$ to this variety is of the form $z\m\k(x,z)+c$ where $c$ is a constant. It suffices to show that the linear
function $z\m\k(x,z)$ on $\fn_{\z\i}$ is not identically zero. If it were identically zero, then $\k(x,\fn)=0$ 
hence $x\in\fp$ and $x\in\fp_\z$ contradicting our assumptions. The lemma is proved.

\subhead A.7\endsubhead 
Assume now that $B$ is an orbital complex (on $\fg_\z$). We show:

(a) {\it any composition factor of any perverse cohomology sheaf of $\tRes_{\fp_\z}^{\fg_\z}(B)$ is an orbital 
complex on $\fl_\z$.}
\nl
Since each perverse cohomology sheaf of $D:=\tRes_{\fp_\z}^{\fg_\z}(B)$ is $L_K$-equivariant, it is enough to show
that $\supp(D)$ is contained in the union of finitely many $L_K$-orbits on $\fl_\z$. Since the intersection of any
$L$-orbit on $\fl$ with $\fl_\z$ is a union of finitely many $L_K$-orbits, it is enough to show that $\supp(D)$ is
contained in the union of finitely many $L$-orbits in $\fl$. Clearly, $\supp(D)$ is contained in 
$\p(\supp(B)\cap\fp)$ hence it is enough to show that $\p(\supp(B)\cap\fp)$ is contained in the union of finitely
many $L$-orbits in $\fl$. Since $\supp(B)$ is the union of finitely many $K$-orbits, it is contained in the union
of finitely many $G$-orbits on $\fg$, and it is enough to note that for any $G$-orbit $\co$ in $\fg$, the set
$\p(\co\cap\fp)$ is the union of finitely many $L$-orbits in $\fl$. This completes the proof of (a).

Next we assume that $B$ is an antiorbital complex (on $\fg_\z$). We show:

(b) {\it any composition factor of any perverse cohomology sheaf of $\tRes_{\fp_\z}^{\fg_\z}(B)$ is an 
antiorbital complex on $\fl_\z$.}
\nl
Let $R$ be a composition factor of a perverse cohomology sheaf of $\tRes_{\fp_\z}^{\fg_\z}(B)$. By assumption we
have $B=\fF(B_1)$ where $B_1$ is an orbital complex on $\fg_{\z\i}$. Using A.5 with $B$ replaced by $B_1$ and 
$\z$ by $\z\i$ we see that $\tRes_{\fp_\z}^{\fg_\z}(B)=\fF(\tRes_{\fp_{\z\i}}^{\fg_{\z\i}}(B_1))$. Let $R_1$ be a
simple perverse sheaf on $\fl_{\z\i}$ such that $R=\fF(R_1)$. Then $R_1$ is a composition factor of a perverse 
cohomology sheaf of $\tRes_{\fp_{\z\i}}^{\fg_{\z\i}}(B_1)$. Using (a) with $B$ replaced by $B_1$ and $\z$ by 
$\z\i$ we see that $R_1$ is orbital. Hence $R$ is antiorbital and (b) is proved.

Finally we assume that $B$ is a biorbital complex (on $\fg_\z$). Combining (a),(b) we obtain:

(c) {\it any composition factor of any perverse cohomology sheaf of $\tRes_{\fp_\z}^{\fg_\z}(B)$ is a biorbital
complex on $\fl_\z$.}

\enddocument